\RequirePackage{fix-cm}
\documentclass[smallcondensed]{svjour3}     
\smartqed  
\usepackage{graphicx}
\usepackage{amsmath}
\usepackage{amssymb}
\usepackage{verbatim}
\usepackage{color}
\usepackage{enumerate}
\usepackage{bbold}
\usepackage[normalem]{ulem}

\newcommand{\ord}[1]{\mathcal{O}\left({#1}\right)}
\renewcommand{\rho}{\varrho}
\newcommand{\im}{\mathbb{i}}
\newcommand{\fvec}[1]{\mathbf{#1}}
\newcommand{\id}{\mathbb{1}} 
\newcommand{\del}{\partial}
\newcommand{\veccc}[3]{\left ( \begin{array}{c}#1\\#2\\#3\\ \end{array}\right )}
\renewcommand{\div}{\nabla \cdot}
\newcommand{\Reyn}{\text{Re}}
\newcommand{\aap}{Astronomy \& Astrophysics\ }

\begin{document}

\title{A numerical scheme for the compressible low-Mach number regime
  of ideal fluid dynamics}

\author{Wasilij Barsukow                    \and
            Philipp V.\ F.\ Edelmann         \and
            Christian Klingenberg              \and
            Fabian Miczek                          \and
            Friedrich K.\ R{\"o}pke             
}

\institute{W.\ Barsukow, C.\ Klingenberg, F.\ Miczek \at
	      W\"urzburg University\\
              Institute for Mathematics \\
              Emil-Fischer-Stra\ss e 40\\
              97074 W\"urzburg\\
              Germany\\
              \email{w.barsukow@mathematik.uni-wuerzburg.de}           
           \and
           Ph.\ V.\ F.\ Edelmann, F.\ K.\ R{\"o}pke \at
              Heidelberg Institute for Theoretical Studies\\
              Schloss-Wolfsbrunnenweg 35\\
              69118 Heidelberg\\
              Germany
}

\date{Received: date / Accepted: date}

\maketitle

\begin{abstract}
Based on the Roe solver a new technique that
allows to correctly re\-present low Mach number flows with a
discretization of the compressible Euler equations was proposed in \cite{miczek2015a}. We analyze properties of
this scheme and demonstrate that its limit yields a discretization of the continuous limit system. Furthermore we perform a linear stability analysis for the case of explicit time integration and study the performance of the scheme under implicit time integration via the evolution of its condition number. A numerical implementation demonstrates
the capabilities of the scheme on the example of the Gresho vortex which
can be accurately followed down to Mach numbers of $\sim 10^{-10}$.
\end{abstract}

\section{Introduction}
\label{intro}

This paper concerns itself with finite volume schemes for the compressible Euler equations, in regimes where the Mach number may become both high and quite low. When the Mach number is of order one, modern shock capturing methods are able to resolve discontinuities and other complex structures with high numerical resolution. For the Mach number going to zero the solutions to the compressible Euler equations
\begin{eqnarray}
 \del_t \rho + \nabla \cdot (\rho \fvec v) &= 0 \label{eq:euler1}\\
 \del_t (\rho \fvec v) + \nabla \cdot \left(\rho\fvec v \otimes \fvec
  v + p\cdot \id\right ) &= 0 \label{eq:euler2}\\
 \del_t E + \nabla \cdot [\fvec v(E+p)] &= 0 \label{eq:euler3}
\end{eqnarray}
tend to the solutions of the incompressible Euler equations. This has been demonstrated by e.g. \cite{klainerman1981a,metivier2001a}. They found that the various functions (pressure, density, etc.) converge to those encountered in the incompressible setting at different rates in the Mach number. A numerical method needs to take account of this.

For the compressible Euler equations the CFL stability criterion of an explicit time discretization requires the time step to be very small for small Mach numbers. This is due to the fact that in this regime the sound waves are much faster than the advection of the flow. Additionally to this stiffness in time, shock-capturing methods show an excessive diffusion which completely deteriorates the solution when the methods are applied to flows with low Mach numbers.

In order to deal with these problems in the literature one finds two approaches:				

\begin{itemize}
\item In one approach the flux function of the finite volume method is modified. The idea is to adapt the flux to the low Mach situation. Recent suggestions of flux modification include \cite{dellacherie2010a,rieper2011a,osswald2015a,chalons2016a}. This is a method that was initially proposed by Eli Turkel (\cite{turkel1999a}) for the calculation of steady state flows and was subsequently extended. These methods however are not well suited for flows where low-Mach flows occur simultaneously with flows that have speeds comparable to the sound speed.
\item In another approach Klein (\cite{klein1995a}) devised an algorithm that keeps track of different orders in the asymptotic expansion of the pressure. The idea is to split the system into two parts. One of them involves a slow, nonlinear and conservative hyperbolic system adequate for the use of modern shock capturing methods mentioned above, and the other is a linear hyperbolic system which contains the stiff acoustic dynamics, which is to be solved implicitly. Recent developments for all-speed schemes of this sort are
\cite{cordier2012a,haack2012a,degond2011a}.
In summary, this leads to a hybrid scheme, partly implicit in time, partly explicit. 
\end{itemize}

In this paper we are inspired by the first approach. It is based on a recently published astrophysical paper \cite{miczek2015a}, where the authors
propose a new hydrodynamics solver based on modifying the diffusion matrix of the Roe scheme. In spirit this is not dissimilar to \cite{turkel1999a}, where these modifications were referred to as flux-preconditioning for historical reasons. In more extreme astrophysical situations, however, the schemes proposed there may fail. \cite{miczek2015a} demonstrated that the flux
function resulting from previous preconditioning techniques may show inconsistencies in certain applications. With their new scheme
a consistent scaling is achieved.

Since this new technique may be useful for applications outside the
astrophysical context, in this paper we follow up on the approach taken by
\cite{miczek2015a} and analyze its properties in detail. The requirements for an all Mach number finite volume scheme, inferred from the limit of the continuous system, are

\begin{enumerate}[(i)]
\item the numerical method for the compressible Euler equations converges formally to a discretization of the incompressible equations \label{req:1}
\item the numerical evolution of the kinetic energy near the incompressible regime is independent of the Mach number\label{req:2}
\end{enumerate}
In addition for an efficient numerical method we require
\begin{enumerate}[(i)]
\setcounter{enumi}{2}
\item linear stability of the scheme when subject to explicit time integration\label{req:3}
\item efficient implicit time integration \label{req:4}
\end{enumerate}

The first two requirements are inspired by the properties of the limit at continuous level which will be formulated in Sect. \ref{sect:lm}. 
After introducing a discretization of the Euler equations these requirements will be given a shape that is reasonable for numerical applications. In Sect.~\ref{sect:precond} we give a more extended motivation for the form of the proposed modification of the flux function. For this method, the requirement (\ref{req:1}) from above is investigated in Sect.~\ref{sect:ap} by pursuing the question of whether the
technique qualifies as an asymptotic-preserving scheme and whether additionally properties beyond a consistent discretization of the limit equations are needed. With numerical experiments we demonstrate that our scheme 
yields satisfactory results for flows down to very low Mach numbers, thus giving evidence that it
complies with the above requirement (\ref{req:1}) and (\ref{req:2}). Because we only modify the diffusion matrix
affecting solely the spatial discretization of the equation, we initially
employ the method of lines. For practical implementations, this raises the question of an
appropriate strategy for time discretization. \cite{miczek2015a} applied their method in both explicit and
implicit time discretization. We discuss stability of the scheme in
explicit time discretization in Sect.~\ref{sect:stabil}, consistent with requirement (\ref{req:3}), and comment on its efficiency in implicit time discretization, which allows to cover extended periods of time, in Sect.~\ref{sect:implicit} (requirement (\ref{req:4})).

\section{Fluid dynamics in the low Mach number limit}
\label{sect:lm}

The solutions to
\begin{eqnarray}
 \del_t \rho + \nabla \cdot (\rho {\fvec v} ) &= 0 \label{eq:eulerrescaled1}\\
 \del_t (\rho \fvec v) + \nabla \cdot \left(\rho {\fvec v} \otimes \fvec v + \frac{p}{M^2}\cdot \id\right ) &= 0 \label{eq:eulerrescaled2}\\
 \del_t E + \nabla \cdot \left [\fvec v(E+p)\right ] &= 0 \label{eq:eulerrescaled3}
\end{eqnarray}
tend to solutions of the incompressible Euler equations as $M \in \mathbb R^+$ tends to zero, i.e. in the limit of low Mach numbers
\cite{ebin1977a,klainerman1981a,ukai1986a,asano1987a,isozaki1987a,kreiss1991a,schochet1994a,metivier2001a}. This can formally be seen by expanding all quantities as series in $M$, e.g. for the pressure this would give
\begin{equation}
 p(x, t) = p^{(0)}(x, t) + Mp^{(1)}(x, t) + M^2 p^{(2)}(x, t) + \mathcal{O}(M^3).
\end{equation}
Inserting these into the above equations, collecting order by order and assuming impermeable boundaries gives
\begin{align}
 p^{(0)} &= \mathrm{const} \label{eq:constp0},\\
 p^{(1)} &= \mathrm{const} \label{eq:constp1},\\
 (\nabla \cdot \fvec v)^{(0)} &= 0, \label{eq:zerodiv}
\end{align}
and
\begin{eqnarray}
 \del_t \rho^{(0)} + \fvec v^{(0)} \cdot \nabla \rho^{(0)} &= 0,\\
 \del_t \fvec v^{(0)} + (\fvec v^{(0)} \cdot \nabla) \fvec v^{(0)}  + \nabla p^{(2)}/\rho^{(0)} &= 0.
\end{eqnarray}
These equations describe incompressible flows. Conditions \eqref{eq:constp0}, \eqref{eq:constp1} and \eqref{eq:zerodiv} are true for every time. Initial data that fulfill them are called {\sl well-prepared}. Not well-prepared initial data may lead to an incompressible flow as well, but then an initial disturbance is produced.

The equation for the kinetic energy $E_\text{kin} = \frac12 \rho |\fvec v|^2$ can be rewritten as
\begin{align}
 \del_t E_\text{kin} + \div \left[ \fvec v \left( E_\text{kin} + \frac{p}{M^2} \right ) \right ] &= \frac{p}{M^2} \div \fvec v.
\end{align}
The source term vanishes for incompressible flows and in this case the kinetic energy becomes a conserved quantity. For compressible flows, this is true in the limit $M\to0$ as well, despite of $\frac{\div v}{M^2} \not\in \mathcal{O}(M)$. Expanding the quantities and using \eqref{eq:constp0} and \eqref{eq:constp1} makes the terms proportional to $\frac{1}{M}$ or $\frac{1}{M^2}$ cancel and gives
\begin{align}
 \del_t E_\text{kin} + \div \left[ \fvec v \left( E_\text{kin} + p^{(2)}  \right ) \right ]  &= p^{(2)} \div \fvec v + \mathcal{O}(M).
\end{align}
Now the source term indeed is $\mathcal{O}(M)$ and the kinetic energy can be seen to become a conserved quantity in the limit $M\to0$.

We are thus led to phrase the requirements (\ref{req:1}) and (\ref{req:2}) from the Introduction in a way as they hold at continuous level:
\begin{enumerate}[(i)]
\item If the initial data for the compressible, homogeneous Euler equations are chosen to have spatial pressure fluctuations scale with $\mathcal{O}(M^2)$ and the divergence of the velocity field scale with $\mathcal{O}(M)$, then the solution converges to the solution of the incompressible Euler equations in the limit $M \to 0$, with only these pressure fluctuations playing the role of the dynamic pressure.
\item For solutions to the compressible, homogeneous Euler equations in the low Mach number limit, the total kinetic energy is conserved.
\end{enumerate}
They will be stated again in Section \ref{sect:precond} in a form adapted to the discrete system. They will thus become requirements for a numerical scheme to be able to capture flows in the low Mach number regime.

\section{Spatial discretization and modification of the diffusion matrix}
\label{sect:precond}

\subsection{Finite volume schemes for conservation laws in the low Mach limit}
Consider a system of conservation laws
\begin{align}
 \del_t \fvec U + \del_x \vec F^{(x)}(\fvec U) + \del_y \vec F^{(y)}(\fvec U) + \del_z \vec F^{(z)}(\fvec U) &= 0 
\end{align}
with $\fvec U$ being the vector of conserved quantities ($\fvec U = (\rho, \rho \vec v, E)^\text T$ in the case of the Euler equations) and $\vec F^{(x)}$, $\vec F^{(y)}$, $\vec F^{(z)}$ the flux in $x$-, $y$- and $z$-direction, respectively.

Its numerical solution of the compressible Euler equations
\eqref{eq:euler1}--\eqref{eq:euler3} is achieved by applying the Godunov method on a Cartesian computational grid
with a uniform spacing $\Delta x$, $\Delta y$, $\Delta z$. The conserved quantities stored in a cell $(i,j,k)$ are denoted by $\fvec{U}_{i,j,k}$, and a numerical flux through the interface between cells $(i, j, k)$ and $(i+1, j, k)$ by $\fvec{F}^{(x)}_{i+1/2,j,k}$. Written as a semi-discrete system the finite volume scheme amounts to
\begin{align}
  \label{eq:finitevolume-disc}
  \frac{\partial}{\partial t} \fvec{U}_{i,j,k} + 
\frac{1}{\Delta x} &\left(  \fvec{F}^{(x)}_{i+1/2,j,k} - \fvec{F}^{(x)}_{i-1/2,j,k} \right) + \\
      \frac{1}{\Delta y} &\left(\fvec{F}^{(y)}_{i,j+1/2,k} - \fvec{F}^{(y)}_{i,j-1/2,k} \right) + \\
      \frac{1}{\Delta z} &\left(\fvec{F}^{(z)}_{i,j,k+1/2} - \fvec{F}^{(z)}_{i,j,k-1/2} \right)     
    = 0.
\end{align}

It was suggested by Roe \cite{roe1981a} to choose the numerical flux as
\begin{align}
  \label{eq:roeflux-unpre}
  \fvec{F}_{l+1/2} = \frac{1}{2} & \left[
    \fvec{F}\left( \fvec{U}^\text L_{l+1/2} \right) +
    \fvec{F}\left( \fvec{U}^\text R_{l+1/2} \right)  -
    \left| A_{\text{Roe}} \right| \left(
      \fvec{U}^\text R_{l+1/2} - \fvec{U}^\text L_{l+1/2}
    \right)
  \right],
\end{align}
with the matrix $A_{\text{Roe}}$ resulting from the solution of a linearized Riemann problem at the interface of the
computational cells. It is the absolute value of the Jacobian~$A$, performed on the eigenvalues, and evaluated in the Roe-average state. This term
ensures upwinding and introduces an artificial viscosity that
stabilizes the scheme. The indices $\text L$ and $\text R$ denote the states to the left and to the right of the cell
interface as determined from an appropriate reconstruction procedure. In multi-dimensional context this flux is applied direction by direction. One may introduce other matrices instead of $|A_\text{Roe}|$, and we will refer to them as diffusion, or upwind artificial viscosity matrices.

As argued in Ref.~\cite{miczek2015a}, the problem arising for this
approach in the low Mach number limit is that the upwind artificial
viscosity dominates all the other terms in the limit of small Mach numbers even if the initial data are well-prepared. In particular it excites spatial pressure fluctuations $\mathcal O(M)$. One thus is led to require for a numerical scheme able to maintain low Mach number flows:
\label{sec:req-num}
\begin{enumerate}[(i)]
\item Considering the limit $M \to 0$, and having initial data for the compressible, homogeneous Euler equations chosen to have spatial pressure fluctuations scale with $\mathcal{O}(M^2)$, then this shall hold for the data at late times as well. 
\label{req:num1}
\item Numerical solutions to the compressible, homogeneous Euler equations in the low Mach number limit and with fixed discretization shall display a (numerical) dissipation of kinetic energy that is $\mathcal O(1)$ in the limit $M \to 0$. In particular this means that the dissipation shall not grow with decreasing $M$. \label{req:num2}
\end{enumerate}
These are modified versions of the findings mentioned in Sect. \ref{sect:lm}, reinterpreted in regard to numerical methods. The emphasis on a fixed discretization is due to the ubiquitous observation that with usual finite volume methods for fixed $M$ the dissipation is effectively reduced by increasing the spatial resolution. However this is neither efficient with respect to the invested computation time nor does it touch the root of the problem, the artefacts just reappearing on finer scales again. 

For numerical methods we require additionally
\begin{enumerate}[(i)]
 \setcounter{enumi}{2}
 \item linear stability under explicit time integration \label{req:num3}
 \item efficient implementation of implicit time integration \label{req:num4}
\end{enumerate}
After the discussion of the proposed scheme in this Section, we address (\ref{req:num1}) and (\ref{req:num2}) in Section \ref{sect:ap}, (\ref{req:num3}) in Section \ref{sect:stabil} and (\ref{req:num4}) in Section \ref{sect:implicit}.

The central flux,
\begin{equation} \label{eq:central_flux}
\vec{F}_{l+1/2} =  \frac{1}{2}\left[ \vec{F}(\vec{U}^L_{l+1/2}) + \vec{F}(\vec{U}^R_{l+1/2})\right].
\end{equation}
would be a choice complying with the low Mach number limit, but it lacks
stability in explicit time discretization.

Numerical tests, as
shown in Figure~3 of \cite{miczek2015a}, suggest that it is stable in
the implicit case. A small growth in kinetic energy is observed, which
is, while not being in contradiction with the basic conservation laws, in
contradiction with thermodynamics and therefore less suited to many
practical applications.

Replacing $\left| A_{\text{Roe}} \right|$ by the modified diffusion matrix
\begin{align}
  P^{-1}\left|P A \right| \label{eq:roeflux-pre}
\end{align}
with $P$ a suitable invertible matrix, and $A$ the Jacobian,
is observed (\cite{turkel1999a}) to improve the numerical solutions in
the low Mach number regime. Miczek et al.\ \cite{miczek2015a} argue that this
is because $P$ can be chosen to correct the scaling behavior with low Mach
number found in the diffusion matrix. For historical reasons $P$ is called preconditioning matrix; we will refer to it as the modifying matrix. A widely used one is due to Weiss and Smith
\cite{weiss1995a}, in primitive variables 
\begin{equation}
\label{eq:pc_weiss}
P_\text{prim} =
  \begin{pmatrix}
    1 & 0 & 0 & 0 & \frac{\mu^2 - 1}{c^2} \\
    0 & 1 & 0 & 0 & 0 \\
    0 & 0 & 1 & 0 & 0 \\
    0 & 0 & 0 & 1 & 0 \\
    0 & 0 & 0 & 0 & \mu^2 
  \end{pmatrix},
\end{equation}
with the parameter \begin{align}\mu = \min[1,\max(M_\text{loc},M_\text{cut})],\label{eq:mudefinition}\end{align} which should scale
with the local Mach number $M_\text{loc}$. Its lower limit,
$M_\text{cut}$, avoids singularity of the matrix.

It corrects the scaling behavior of almost all entries in the diffusion matrix. Indeed, when used with the homogeneous Euler equations it shows very low, Mach number independent dissipation in the low Mach limit. However, when integrated implicitly in time, this scheme displays unsatisfactory behavior of the condition number as discussed later.

\subsection{Low Mach modifications in presence of gravity source terms}

A problem with the particular choice \eqref{eq:pc_weiss}
of the dissipation matrix arises if it is used in the presence of certain
source terms, e.g. gravity.
Here the lowest order in the expansion of the pressure in powers of $M$
is not constant. The diffusion matrix obtained when using
\eqref{eq:pc_weiss} has an entry $\mathcal{O}(1/M^2)$ in the energy row.
Therefore with a spatially varying background this introduces strong
diffusion. An example in which such configurations arise are hydrostatic equilibria in presence of gravity. Here the rescaled Euler equations, \eqref{eq:eulerrescaled1}--\eqref{eq:eulerrescaled3}, change to
\newcommand{\fr}{\mathit{Fr}}
\begin{align}
 \del_t \rho + \nabla \cdot (\rho \fvec v) &= 0 \label{eq:euler1-g},\\
 \del_t (\rho \fvec v) + \nabla \cdot \left(\rho\fvec v \otimes \fvec
 v + \frac{p}{M^2}\cdot \id\right ) &= \frac{\rho}{\fr^2} \fvec{g} \label{eq:euler2-g},\\
 \del_t E + \nabla \cdot [\fvec v(E+p)] &= \frac{M^2}{\fr ^2}\rho \fvec{g} \cdot \fvec{v}. \label{eq:euler3-g}
\end{align}
The nondimensional Froude number $\fr $ appears in presence of gravity source terms. For the purposes of this article
we only consider the case of $\fr =M$.

This new system admits static solutions (stationary and with $\fvec v = 0$) called
\emph{hydrostatic equilibria} which are governed by the condition $\nabla p = \rho \fvec g$. 

To point out and compare the properties of the numerical schemes in presence of this source term, consider a special but important solution in one spatial dimension. In
the case of a constant temperature $T$, spatially and temporally
constant $\fvec g = -g \fvec e_x$ pointing into the negative
$x$-direction and the ideal gas equation of state, hydrostatic equilibrium has the form
\begin{equation}
  \rho(x) = \rho_0 \exp\left(-\frac{gx}{T}\right), \quad p(x) = \rho(x)T, \quad \rho_0 := \rho(0)
\end{equation}
We can compute the discrete fluxes using
\eqref{eq:roeflux-pre} on a uniform grid with spacing $\Delta x$
using constant reconstruction. The change in grid cell~$i$ is given by
\begin{equation}
  \label{eq:spatial-residual}
  -\del_t \fvec U_i = \begin{pmatrix}0\\
  \frac{p_{i+1} - p_{i-1}}{2 \Delta x M^2} + \rho_i
  g\\
  0
  \end{pmatrix}
  + \frac{1}{2\Delta x} \left(-\fvec D_{i+1/2} + \fvec D_{i-1/2}\right).
\end{equation}
with $\fvec D_{i+1/2} = (P^{-1}|PA)_{\text{Roe},i+1/2} \; (\fvec{U}_{i+1} - \fvec{U}_{i})$.

For the first expression, which corresponds to the central flux,
Eq.~\eqref{eq:central_flux}, we get
\begin{align}
  \label{eq:hystat-p-diff}
  \frac{p_{i+1}-p_{i-1}}{2\Delta x } &= \frac{\rho_0 
  T}{2\Delta x } \left[\exp\left(-\frac{g(i+1)\Delta x}{T}\right) -
  \exp\left(-\frac{g(i-1)\Delta x}{T}\right)\right]\\
  \nonumber
  &=- \rho_i g \left[1 + \frac{g^2 \Delta x^2}{6 T^2} +
  \ord{\Delta x^4}\right].
\end{align}
which cancels with the cell-centered discretization of gravity up to
order of $\Delta x^2$.

The effect of the numerical dissipation
term~$\fvec D_{i\pm 1/2}$ can
most conveniently be evaluated in primitive variables ($\rho, \vec v, p$). For the modifying matrix $P$ in Eq.~\eqref{eq:pc_weiss} the
contribution is
\begin{equation}
  \label{eq:hystat-diff-turkel}
\frac{\fvec D_{i+1/2}}{2\Delta x} = \begin{pmatrix}
    1\\
    0\\
    c^2
  \end{pmatrix}
  \frac{\Delta p}{2 \Delta x c M_\text{cut}},
\end{equation}
with $\Delta p = p_{i+1} - p_i$ and the local speed of sound~$c=\sqrt{\gamma p / \rho}$ with the
$\gamma$ from Eq.~\eqref{eq:eos}. As we consider hydrostatic solutions here, the bounded local Mach number $\mu$ \eqref{eq:mudefinition} was set to its lower limit $M_\text{cut}$. Recall that $M_\text{cut}$ is normally
chosen to a value well below the Mach number in the considered flow
field, just large enough to avoid infinite values when dividing by $\mu$ in zero-velocity
regions, while maintaining the positive effect of the modification of the diffusion matrix at
reasonably low Mach numbers. Eq. \eqref{eq:hystat-diff-turkel} reveals a fundamental problem
when using the matrix \eqref{eq:pc_weiss} -- its contribution becomes extremely
large in regions with very low Mach numbers if a pressure gradient is present.
While this is generally not the case in low Mach number flows when the homogeneous Euler equations are solved, it can
happen, if gravity or other source terms are involved.

Because of these problems, Miczek et al.\ \cite{miczek2015a} suggested
a new modifying matrix~$P$. In entropy variables it takes the form,
  \begin{equation}
    \label{eq:pc_miczek_sym}
    P_\text{entr} =
    \begin{pmatrix}
      1          & n_x \delta & n_y \delta & n_z \delta & 0 \\[1ex]
      -n_x \delta & 1          & 0          & 0          & 0 \\[1ex]
      -n_y \delta & 0          & 1          & 0          & 0 \\[1ex]
      -n_z \delta & 0          & 0          & 1          & 0 \\[1ex]
      0          & 0          & 0          & 0          & 1 
    \end{pmatrix},
  \end{equation}
  with $\delta = \frac{1}{\min(1, \max(M_\text{loc},M_\text{cut}))} - 1$.
In primitive variables it is
\begin{equation}
  \label{eq:pc_miczek_prim}
  P_\text{prim} =
  \begin{pmatrix}
      1 &
      n_x \frac{\rho \delta M}{c} &
      n_y \frac{\rho \delta M}{c} &
      n_z \frac{\rho \delta M}{c} & 0 \\[1ex]
      0 & 1 & 0 & 0 & -n_x \frac{\delta}{\rho c M}  \\[1ex]
      0 & 0 & 1 & 0 & -n_y \frac{\delta}{\rho c M}  \\[1ex]
      0 & 0 & 0 & 1 & -n_z \frac{\delta}{\rho c M}  \\[1ex]
      0 &
      n_x \rho c \delta M &
      n_y \rho c \delta M &
      n_z \rho c \delta M & 1
  \end{pmatrix}.
\end{equation}
As for the modifying matrix from Eq.~\eqref{eq:pc_weiss}, the
definition of $\delta$ ensures that the scheme reverts back to the
original Roe scheme when the local Mach number reaches 1. Miczek et al.\
\cite{miczek2015a} show that the scaling with $M$ of the diffusion
matrix of this scheme is fully consistent with the flux Jacobian. The flux at any cell edge is given, similarly to \eqref{eq:roeflux-unpre} by
\begin{align*}
  \fvec{F}_{l+1/2} = \frac{1}{2} & \left[
    \fvec{F}\left( \fvec{U}^\text L_{l+1/2} \right) +
    \fvec{F}\left( \fvec{U}^\text R_{l+1/2} \right)  -
    P^{-1}\left|P A \right| \left(
      \fvec{U}^\text R_{l+1/2} - \fvec{U}^\text L_{l+1/2}
    \right)
  \right]
\end{align*}

We can perform the analysis of the discretized fluxes in hydrostatic
equilibrium for this solver, too. Equation \eqref{eq:hystat-p-diff} is
identical for both. Instead of the dissipation term from
Eq.~\eqref{eq:hystat-diff-turkel}, we get in primitive variables
\begin{equation}
  \frac{\fvec D_{i+1/2}}{2\Delta x} = \begin{pmatrix}
    M_\text{cut} c^{-1}\\
    \frac{M_\text{cut}-1}{\rho M}\\
    c M_\text{cut}
  \end{pmatrix}
  \frac{\Delta p}{2 \Delta x M \sqrt{1 - 2 M_\text{cut} + 2
  M_\text{cut}^2 }}.
\end{equation}
This expression overcomes the problems of the preconditioner in Eq.~\eqref{eq:pc_weiss} because its dissipation term does not grow when lowering $M_\text{cut}$.
Therefore the new method avoids the problems encountered for previously suggested modification matrices in presence of gravity.

Solutions to the Euler equations augmented by a gravity source term need not in general be or converge to a static equilibrium. However there are indeed a lot of interesting applications related to such equilibria. The ability of a scheme to preserve them up to machine precision is a challenging additional requirement, an implementation of which however is not the topic of this paper but subject of ongoing work.

\section{Asymptotic behavior of the numerical method}
\label{sect:ap}

\subsection{Asymptotic analysis of the semi-discrete scheme}
It has been demonstrated in \cite{miczek2015a} that the new
way of modifying the diffusion matrix (Eq. \eqref{eq:pc_miczek_sym} or \eqref{eq:pc_miczek_prim}) ensures that the diffusive part $$P^{-1}\left|P A \right| \left(
      \fvec{U}^\text R_{l+1/2} - \fvec{U}^\text L_{l+1/2}
    \right)$$ does not dominate the numerical flux function at
low Mach numbers. In the context of asymptotic preserving schemes it has been found useful to analyze the limit of the discrete system, in a way analogous to what has been done in Sect.~\ref{sect:lm} in the continuous case. With, for simplicity, a piecewise constant reconstruction, the numerical flux in $x$-direction is given by
\newcommand{\blue}{}
\newcommand{\num}{}
\newcommand{\veccccc}[5]{\left ( \begin{array}{c}#1\\#2\\#3\\#4\\#5 \end{array}\right )}
\renewcommand{\boxed}{}
\begin{equation}
 \fvec F_{i+\frac12} = \frac{1}{2}\left [\fvec F(\fvec U_{i+1}) + \fvec F(\fvec U_{i}) \right ] - \frac{1}{2} P^{-1} |P A| (\fvec U_{i+1} - \fvec U_{i} ), \label{eq:preconditionedflux}
\end{equation}
where $P^{-1} |P A|$ is evaluated in the Roe state, and $A$ is the Jacobian in $x$-direction (indices for the other directions have been dropped for better readability). The fluxes through the other interfaces can be obtained analogously.

Taking $\delta \in \mathcal{O}\left(\frac{1}{\blue M}\right )$, by
construction the leading order terms in $\blue M$ of the diffusion
matrix are the same as in the Jacobian (\cite{miczek2015a}). In the basis of conserved variables one finds
 \begin{align}
  P^{-1} |P A| &= \frac{1}{\blue M^2} \left( \begin{array}{ccccc} 
  0 & 0   &0  & 0                            & 0 \\ 0 & 0 &0& 0 & \gamma-1 \\ 0 & 0&0 & 0 & 0\\0 & 0&0 & 0 & 0\\0 & 0&0 & 0 & 0\end{array} \right )  + \mathcal{O}(1).
 \end{align}
 A conservative semi-discrete scheme with flux \eqref{eq:preconditionedflux} then is
 \begin{align*}
 0 &=  \del_t \fvec U_{\num i} + \frac{\fvec F_{\num i+\frac12} - \fvec F_{\num i-\frac12}}{\Delta x} + \text{ fluxes through other interfaces} 
 \end{align*}
 and to highest order
 \begin{align*}
 0&= \del_t \fvec U_{\num i} + \frac{1}{2\Delta x} \left [ \frac{1}{\blue M^2} \veccccc{0}{p_{\num i+1} - p_{\num i-1}}{0}{0}{0} - \frac{\gamma-1}{\blue M^2}\veccccc{0}{E_{\num i+1} - 2E_{\num i} + E_{\num i-1}}{0}{0}{0}  \right ]  + \mathcal{O}({\blue M})
 \end{align*}
 
 The two lowest orders can be simplified (for $\ell = 0, 1$ one has $p^{(\ell)} = (\gamma-1) E^{(\ell)}$) to formally yield in the limit $M\to 0$:
 \begin{eqnarray}
  p^{(\ell)}_{\num i} - p^{(\ell)}_{\num i-1} = 0 \qquad \ell = 0,1 \label{eq:pressurediscrete}
 \end{eqnarray}
 
 The rest of the asymptotic analysis is done with the $\mathcal{O}(1)$
 equations, which due to the consistency of the scheme give consistent
 discretizations of the remaining equations in the limit $M\to
 0$. 
 
 Equation \eqref{eq:pressurediscrete} for the Roe solver is, as has been discussed in \cite{guillard1999a,guillard2004a}
 \begin{align}
  \frac{p^{(0)}_{\num i+1} - p^{(0)}_{\num i-1} }{\Delta x} &= 0 \\
  \frac{p^{(1)}_{\num i+1} - p^{(1)}_{\num i-1} }{\Delta x}&= \Delta x \cdot (\text{terms involving $2^\text{nd}$ derivatives of $\rho$, $v$, $e$}) + \mathcal{O}(\Delta x^2).
 \end{align}
 Even though it is also a discretization of $\nabla p^{(\ell)} = 0$, it is not a good approximation for finite values of $\Delta x$ -- contrary to \eqref{eq:pressurediscrete}. As is well-known it is always possible to cure the low Mach number problems by increasing the resolution. This however, as mentioned above, is both impractical and unnecessary. 
 
 In view of the findings we expect our scheme to have pressure perturbations $\mathcal O(M^2)$, consistently with requirement (\ref{req:num1}) from Sect. \ref{sect:precond}.

 \subsection{Numerical results}
 
We demonstrate this result with numerical experiments in which we
simulate a Gresho-vortex setup \cite{gresho1990a}. This is an example of a stationary, incompressible rotating flow around the origin in two spatial dimensions:
\begin{align}
     \vec v &= \vec e_\phi \cdot \begin{cases} 5r & r < 0.2\\ 2 - 5r & r < 0.4 \\ 0 & \text{else} \end{cases} \label{eq:gresho1}\\ 
     p &= \begin{cases} p_\text c + \frac{25}{2} r^2 & r < 0.2 \\ p_\text c + 4\ln(5 r) + 4 - 20 r + \frac{25}{2} r^2 & r < 0.4 \\ p_\text c + 4 \ln 2 - 2 & \text{else}\end{cases}\label{eq:gresho2}
                                                                                         \end{align}
with the uniform density $\rho = 1$ and the pressure in the vortex center $p_\text c = \frac1{\gamma M^2} - \frac12$. Also $r = \sqrt{x^2 + y^2}$ and $\vec e_\phi$ is the azimuthal unit vector in two-dimensional polar coordinates. \\In the compressible setting the flow can be endowed with different maximum Mach numbers by varying the parameter $M$ in the value of the central pressure. Therefore this is an example of a family of solutions, parametrized by a real number $M$, such that $M_\text{loc}$ scales asymptotically as $M$ in the limit $M\to 0$. Here all quantities are understood to be non-rescaled, and one observes for example that $p = \frac{1}{M^2} (\tilde p^{(0)} + M \tilde p^{(1)} + \ldots)$ and so on.

The setup and the employed numerical method are
identical to those presented in Ref.~\cite{miczek2015a}. For the
numerical solution, a fully discretized scheme is necessary and we chose an
implicit time discretization with an advective CFL criterion to
determine the time step size (see Sec. \ref{sect:implicit}). A piecewise linear MUSCL-like
reconstruction without limiters is used
\cite{vanleer1979a,toro2009a}. 

\begin{figure*}
\includegraphics[width=\linewidth]{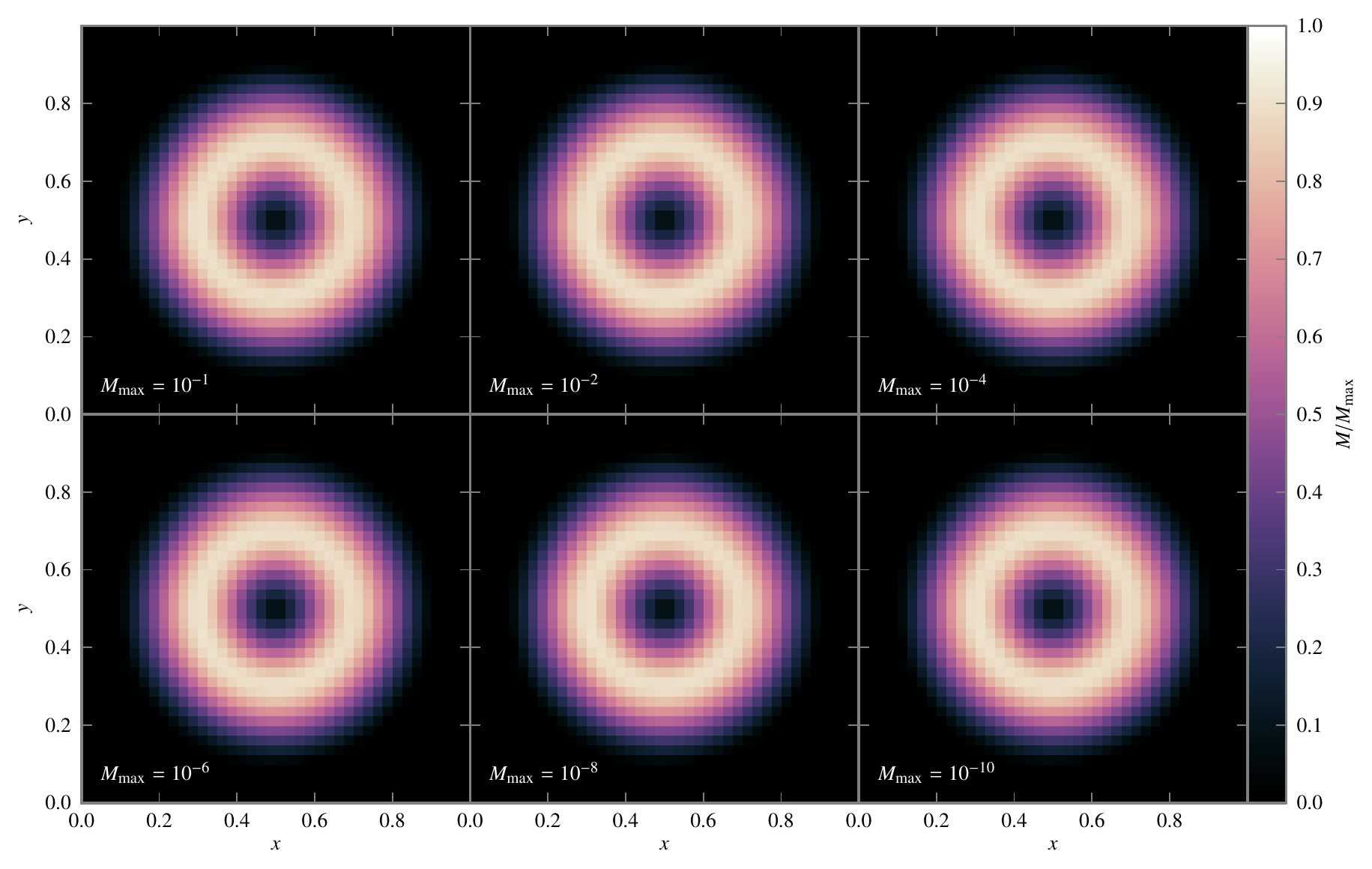}
\caption{Gresho vortex problem advanced over one full revolution with modified
  fluxes for different maximum Mach
  numbers $M_{\mathrm{max}}$ in the setup, as indicated in the plots. Color coded is the
  Mach number relative to the respective $M_{\mathrm{max}}$.}
\label{fig:gresho_multi}
\end{figure*}

As stated above, our goal is to devise
a scheme that represents low Mach number flows well at low numerical
resolution. Therefore, the Gresho vortex is set up on a grid with only
$40 \times 40$ computational cells. We follow the flow over one full revolution of
the vortex and show the results for maximum Mach numbers down to
$10^{-10}$ in Fig.~\ref{fig:gresho_multi}. With the Miczek scheme
\cite{miczek2015a}, the vortex is retained in the simulations even
at the lowest Mach numbers. This contrasts the result obtained with a conventional Roe-type flux function in which
the vortex is significantly blurred after one full revolution at a
maximum Mach number of $10^{-2}$ and completely destroyed for maximum
Mach numbers below $10^{-3}$ as seen in Fig.~\ref{fig:gresho_diff} and \cite{miczek2015a}.

\begin{figure*}
\includegraphics[width=\linewidth]{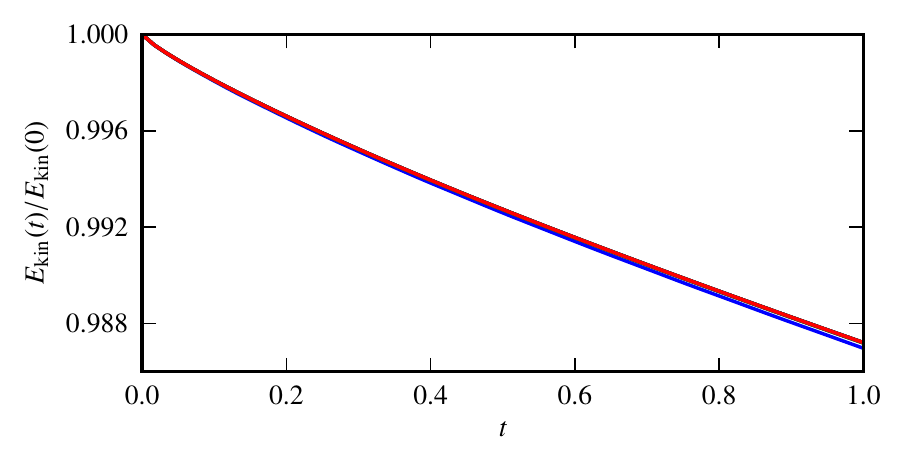}
\caption{Temporal evolution of the total kinetic energy
  $E_{\mathrm{kin}}(t)$ relative to
  its initial value $E_{\mathrm{kin}}(0)$ in the Gresho vortex problem advanced with
  modified fluxes. The cases for $M_{\mathrm{max}} =
  10^{-n}$, $n = 1, 2, 3, 4, 5, 6, 7, 8, 9, 10$
  are overplotted but indistinguishable.}
\label{fig:gresho_ekin_multi}
\end{figure*}

\begin{figure*}
\includegraphics[width=\linewidth]{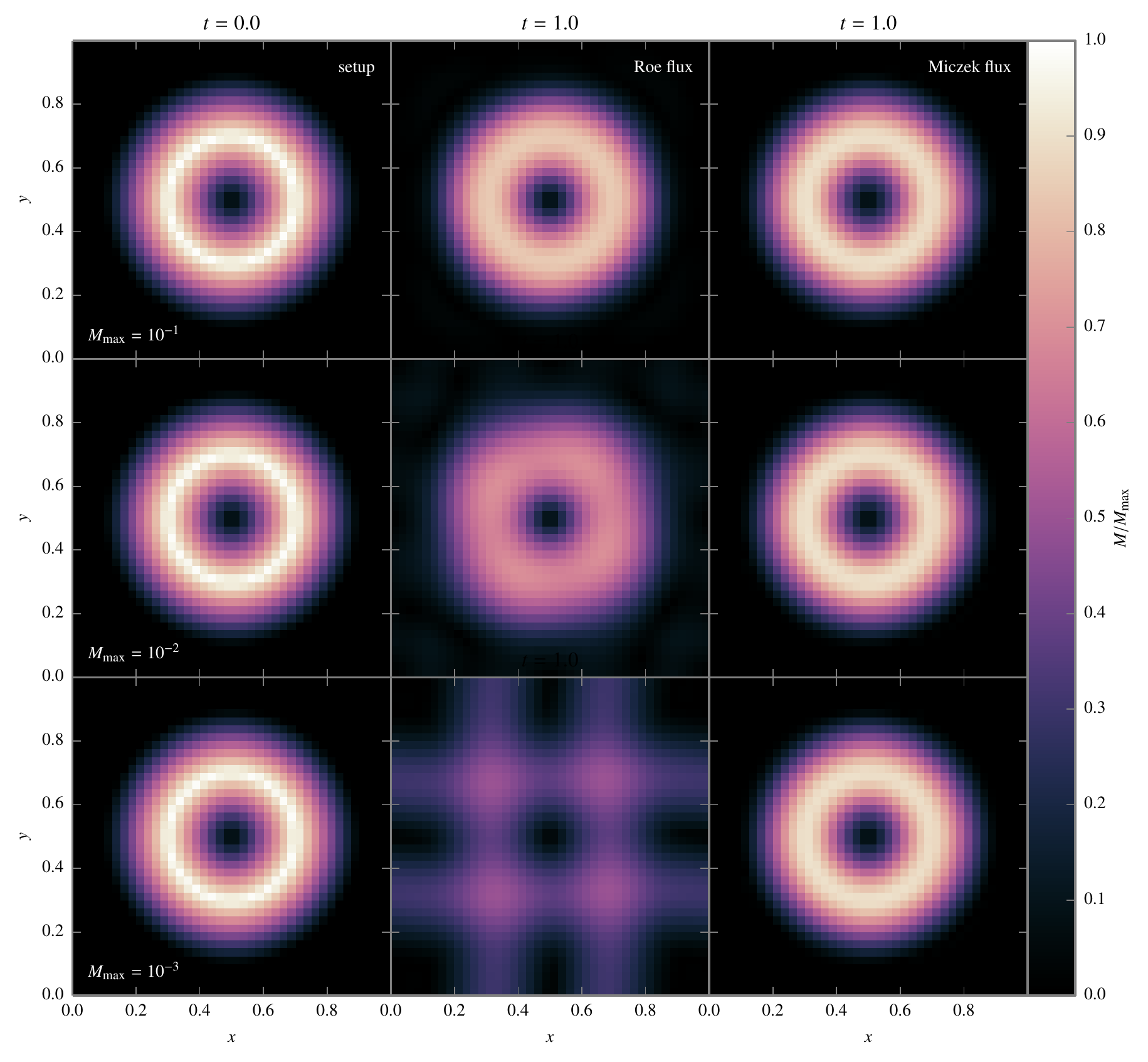}
\caption{Comparison of the Gresho vortex problem at different initial
conditions and computed with different schemes. The rows show the
different initial maximum Mach number. The first column is the initial
condition, the others show the state after one full rotation using 
the unmodified Roe solver,
and the scheme
with the new modifying matrix.
Color coded is the Mach number
normalized to its initial maximum value.}
\label{fig:gresho_diff}
\end{figure*}

The evolution of the total kinetic energy in the simulation domain is
shown in Fig.~\ref{fig:gresho_ekin_multi} and Table~\ref{tab:ekin}. Although the kinetic energy reduces by about 1.3 per cent
over one revolution of the vortex in our setup, this loss is
independent of the Mach number of the flow. This is very much in
contrast to conventional schemes, whose dissipation rate of kinetic
energy increases excessively the lower the Mach numbers get.

At high Mach numbers of about $10^{-1}$ the proposed scheme performs similarly to conventional ones. In view of the results one is however led to the observation that it is at the same time able to reproduce flows at very low Mach
numbers. Its numerical dissipation is not increasing in this limit. We
thus show with our numerical experiments the Miczek et al.\ scheme
\cite{miczek2015a} to fulfill requirements (\ref{req:num1}) and (\ref{req:num2}) as
formulated in Sect.~\ref{sec:req-num}.

In addition to the incompressible flow one still may have sound waves. Since our new scheme is based on a discretization of the full compressible Euler equations \eqref{eq:euler1}--\eqref{eq:euler3} it does not remove them. This is demonstrated in \cite{miczek2015a}, where the example of a sound wave passing through a low Mach vortex is correctly simulated. 

\begin{table}
\caption{Total kinetic energy  $E_{\mathrm{kin}}(1.0)$ after one full
  revolution of the Gresho vortex relative to its initial value
  $E_{\mathrm{kin}}(0.0)$ for different maximum Mach numbers $M_{\mathrm{max}}$.}
\label{tab:ekin}       
\begin{tabular}{rl}
\hline\noalign{\smallskip}
$\log M_{\mathrm{max}}$ & $E_\mathrm{kin}(1.0)/E_\mathrm{kin}(0.0)$\\
\noalign{\smallskip}\hline\noalign{\smallskip}
$-1$  & $0.986974319078$\\  
$-2$  & $0.987185681481$ \\ 
$-3$  & $0.987206395072$ \\ 
$-4$  & $0.987208424676$ \\ 
$-5$  & $0.987208767527$ \\ 
$-6$  & $0.987208721327$ \\ 
$-7$  & $0.987208711049$ \\ 
$-8$  & $0.987208711129$ \\ 
$-9$  & $0.987208710852$ \\ 
$-10$ & $0.987208711987$ \\ 
\noalign{\smallskip}\hline
\end{tabular}
\end{table}

\section{Linear stability of explicit time discretization}
\label{sect:stabil}

\subsection{Stability Analysis}
The correct reproduction of solutions in the low Mach number limit was
achieved by modifying the artificial upwind viscosity matrix -- a term
that was introduced to stabilize the
scheme. This raises the question of the stability of the resulting new
method in explicit time discretization.

The investigation of linear stability with the von Neumann method yields results on the time behavior of Fourier modes for a linear(ized) conservation law. If all of the modes are damped in time, the method is called linearly stable. Surely, a necessary requirement is that the method is stable already in one spatial dimension and when integrated in time by a first order method. For simplicity, the following stability analysis is performed with piecewise constant reconstruction, i.e.\ on a method that is first order in space and time.

Express every quantity $\fvec U^n_i$ by a Fourier series in space ($\im^2 = -1$): \begin{equation} \fvec{U}^n_i = \sum_{k\in\mathbb{Z}} \fvec{U}^n \exp(\im i k \Delta x) \end{equation}  
insert this into the fully discrete scheme ($\nu = \frac{\Delta t}{\Delta x}$)
\begin{equation}
 \fvec U^{n+1}_i = \fvec U^{n}_{i} - \frac{1}{2} \nu \Big[ A (\fvec U^{n}_{i+1} - \fvec U^{n}_{i-1}) - D(\fvec U^{n}_{i+1} - 2\fvec U^{n}_{i} + \fvec U^{n}_{i-1})   \Big ]
\end{equation}
to obtain, by defining $k\Delta x =: \beta$,
\begin{equation}
 \fvec U^{n+1} = \Big \{ \id - \nu \left[ A \im \sin \beta + D(1-\cos \beta)   \right ] \Big \} \fvec U^{n} \label{eq:amplificationmatrix}
\end{equation}
The expression in curly brackets is called amplification matrix. Stability of such iterated linear maps needs all its eigenvalues to be less than 1 in absolute value. In particular, it is considered necessary in \cite{dellacherie2009a} for the absolute value of the acoustic eigenvalues to be strictly less than 1 for the so-called checkerboard-mode $\beta = \pi$. This amounts to non-vanishing eigenvalues of $D$ which is the case for the considered specific choice of $D$, as will be seen later.

Consider the following system
\begin{eqnarray}
 \del_t \veccc{q_1}{q_2}{q_3} + \left(  \begin{array}{ccc} a & a_{12} & 0 \\ 0 &a & a_{23} \\ 0 & a_{32} & a   \end{array} \right ) \del_x \veccc{q_1}{q_2}{q_3} = 0
\end{eqnarray}
which shall be solved with a time-explicit scheme of Roe-type with a diffusion matrix
\begin{equation}
 D = \left(  \begin{array}{ccc} d_{11} & d_{12} & d_{13}\\0 & d_{22} & d_{23} \\ 0 & d_{32} & d_{33} \end{array} \right) 
\end{equation}
The Jacobian of hydrodynamics in one spatial dimension and in primitive variables is of this type. As the stability analysis is linear, the equation may be considered in any variables, not necessarily the conserved ones. Moreover, the study of the eigenvalues $\lambda$ of the amplification matrix in \eqref{eq:amplificationmatrix} equally does not depend on the chosen basis. The property of the diffusion matrix of having just one non-zero entry in the $\rho$-column will be fulfilled by the one appearing here.

The eigenspace decomposes into

\begin{equation}
 1 - \nu ( a \im  \sin \beta +  d_{11}(1-\cos \beta)) = \lambda \label{eq:hydro1dpart}
\end{equation}
and
\begin{eqnarray}
  &[1-\nu ( a \im  \sin \beta +  d_{22}(1-\cos \beta)) - \lambda] [1-\nu ( a \im  \sin \beta +  d_{22}(1-\cos \beta)) - \lambda] \\&= \nu^2 [ a_{32} \im  \sin \beta +  d_{32}(1-\cos \beta)] [  a_{23} \im  \sin \beta +  d_{23}(1-\cos \beta)] \label{eq:hydro2dpart}
\end{eqnarray}
Equation \eqref{eq:hydro1dpart} is easily recognized as a 1-dimensional stability result.
It leads to the stability condition $d_{11} \geq |a|$ and if $d_{11} = |a|$ (as will turn out later in the specific example), then $\nu < \frac{1}{d_{11}}$.

Equation \eqref{eq:hydro2dpart} is just the stability condition for the truncated matrices of a reduced system
\begin{equation}
 A_\mathrm{red}=\left ( \begin{array}{cc} a & a_{23} \\ a_{32} & a \end{array} \right ) \qquad D_\mathrm{red} = \left (\begin{array}{cc} d_{22} & d_{23} \\ d_{32} & d_{33} \end{array} \right )
\end{equation}
Note that the elements $a_{12}$, $d_{12}$ and $d_{13}$ are irrelevant for stability.

Equation \eqref{eq:hydro2dpart} can be rewritten as
\begin{equation}
 \left(1 - \nu \left ( a \im  \sin \beta + \frac{d_{22}+d_{33}}{2} (1-\cos \beta)  \right ) - \lambda \right )^2  = \nu^2 (\mathcal{A}+\mathcal B\im ) 
\end{equation}
with
\begin{align}
 \mathcal A &:= -a_{23} a_{32} \sin^2 \beta + d_{23} d_{32} (1-\cos \beta)^2 + \left( \frac{d_{33} - d_{22} }{2}\right )^2 (1 - \cos \beta)^2\\
 \mathcal B &:= (a_{23} d_{32} + d_{23}  a_{32} )  (1-\cos \beta) \sin \beta\\
 \bar d &:= \frac{d_{22}+d_{33}}{2}
\end{align}
such that
\begin{equation}
 \lambda = 1 - \nu a \im \sin \beta - \nu \bar d (1-\cos \beta)  \pm \nu \sqrt{\mathcal A + \mathcal B \im} 
\end{equation}
The square root is given by
\begin{equation}
 \sqrt{\mathcal A + \mathcal B\im} = \sqrt{\frac{\sqrt{\mathcal A^2 + \mathcal B^2} + \mathcal A}{2}} + \im  \text{\,sgn}(\mathcal B) \cdot \sqrt{\frac{\sqrt{\mathcal A^2 + \mathcal B^2} - \mathcal A}{2}}
\end{equation}
Evaluating $|\lambda|^2 < 1$ leads to
 \begin{align*}
 &\nu < 2 \frac{ \bar d (1-\cos \beta) \mp \sqrt{\frac{\sqrt{\mathcal A^2 + \mathcal B^2} + \mathcal A}{2}} }{\left( \bar d (1-\cos \beta) \mp \sqrt{\frac{\sqrt{\mathcal A^2 + \mathcal B^2} + \mathcal A}{2}} \right )^2 +  \left( a \sin \beta \mp \text{sgn\,}(\mathcal B) \sqrt{\frac{\sqrt{\mathcal A^2 + \mathcal B^2} - \mathcal A}{2}} \right )^2} 
 \end{align*}

 The suggested upwinding matrix from \cite{miczek2015a} is
 \begin{equation}
  \left( \begin{array}{ccc} 
	    |v| & \frac{\rho(-c^2 \delta + c M v  + \delta M^2 v^2)}{c\tau} & -\frac{|v|}{c^2} + \frac{1}{M \tau}\\
	    0 & \frac{c^2}{M \tau} & \frac{c^2 \delta + cMv - \delta M^2 v^2}{cM^2 \rho \tau}\\
	    0 & \frac{c \rho (-c^2 \delta + c M v + \delta M^2 v^2)}{\tau} & \frac{c^2}{M \tau}
	    \end{array}
\right )
 \end{equation}
 with $\tau = \sqrt{c^2 (1 + \delta^2) - \delta^2 M^2 v^2}$.
 
One can investigate the limit of small $M$. For the components of the upwinding matrix one has (having in mind the two cases $\delta \in \mathcal{O}(\frac{1}{M})$ and $\delta \in \mathcal{O}(1)$):
\begin{eqnarray}
 \tau &\sim   c \sqrt{1 + \delta^2} 
 \end{eqnarray}
 \begin{eqnarray}
 \bar d \sim \frac{c}{\sqrt{1 + \delta^2} M} \quad 
 d_{12} \sim \frac{\delta }{\sqrt{1 + \delta^2} M^2} \quad 
 d_{21} \sim -\frac{c^2 \delta}{\sqrt{1 + \delta^2}}
\end{eqnarray}
Therefore
\begin{align}
 \mathcal A &= - \frac{c^2}{M^2} \left( \sin^2 \beta + \frac{\delta^2  }{ 1 + \delta^2} (1-\cos \beta)^2 \right )\\
 \mathcal B 
 &= \frac{2cv}{  \sqrt{1 + \delta^2} M}   (1-\cos \beta) \sin \beta
\end{align}
where due to a lot of cancellations the exact values were used for $\mathcal B$. Whereas in both cases $\mathcal A \in \mathcal{O}(1/M^2)$, one has
\begin{align}
\mathcal B &\in \mathcal O(1/M) \qquad\text{ if } \delta \in \mathcal{O}(1)    \\ \mathcal B &\in  \mathcal O(1) \qquad\text{ if } \delta \in \mathcal O(1/M) \end{align}
As, $|\mathcal A| = -\mathcal A$, $\sqrt{\mathcal A^2 + \mathcal B^2}-\mathcal A \sim 2|\mathcal A|$ and 
\begin{align}
  \sqrt{\mathcal A^2 + \mathcal B^2}+\mathcal A 
  \sim |\mathcal A| \frac{\mathcal B^2}{2 \mathcal A^2} =  \frac{\mathcal B^2}{2 |\mathcal A|}  \in \begin{cases}  \mathcal{O}(1) & \text{ if } \delta \in \mathcal{O}(1)  \\ \mathcal O(M^2) & \text{ if } \delta \in \mathcal O(1/M) \end{cases}  
\end{align}
The term $\sqrt{\frac{\sqrt{\mathcal A^2 + \mathcal B^2} + \mathcal A}{2}} $ will be compared to 
\begin{align}
\bar d (1-\cos \beta) \in \begin{cases} \mathcal{O}(1/M) &\text{ if } \delta \in \mathcal{O}(1)    \\ \mathcal O(1) &\text{ if } \delta \in \mathcal O(1/M) \end{cases}
\end{align}
and the latter wins in both cases. Therefore
\begin{align}
 \nu &< 2 \frac{ d (1-\cos \beta)  }{d^2 (1-\cos \beta)^2  +  \left( c \sin \beta \mp \text{sgn\,}(B) \sqrt{|A|} \right )^2} \\
 &\sim 2 \frac{ \frac{c}{\sqrt{1 + \delta^2} M} (1-\cos \beta)  }{\frac{c^2}{(1 + \delta^2) M^2} (1-\cos \beta)^2  +   \frac{c^2}{M^2} \left| \sin^2 \beta - \frac{\delta^2  }{ 1 + \delta^2} (1-\cos \beta)^2 \right|} \\
&= \frac{M}{c} \frac{ \frac{2}{\sqrt{1 + \delta^2}}   }{\frac{1}{1 + \delta^2} (1-\cos \beta)  +   \left| (1 + \cos \beta) - \frac{\delta^2  }{ 1 + \delta^2} (1-\cos \beta) \right|} \\
&= \frac{M}{c} \frac{ 2\sqrt{1 + \delta^2}   }{ 1-\cos \beta  +   \left| 1 + (1 + 2\delta^2) \cos \beta \right|}
 \end{align}
Now a minimum over all $\beta \in [0, 2\pi)$ has to be performed in order to obtain the global maximum value of $\nu$. If $\delta \in \mathcal O(1)$ (in particular one might be interested to recover for $\delta = 0$ the usual Roe scheme) then \begin{align}\nu_\text{max} \sim \frac{M}{c}\end{align} 
if a suitable minimizing $\cos \beta_\text{min}$ exists, which is $\mathcal{O}(1)$ (trivially the case for $\delta = 0$).

However if $\delta \in \mathcal{O}(1/M)$, then $|\cos \beta_\text{min}| = 1$ and 
\begin{align}
 \nu_\text{max} &\sim \frac{M}{c} \frac{\sqrt{1 + \delta^2}   }{\delta^2 } \in \mathcal{O}\left(\frac{M^2}{c}\right)
\end{align}

\subsection{Numerical Verification}
\begin{figure*}
  \includegraphics{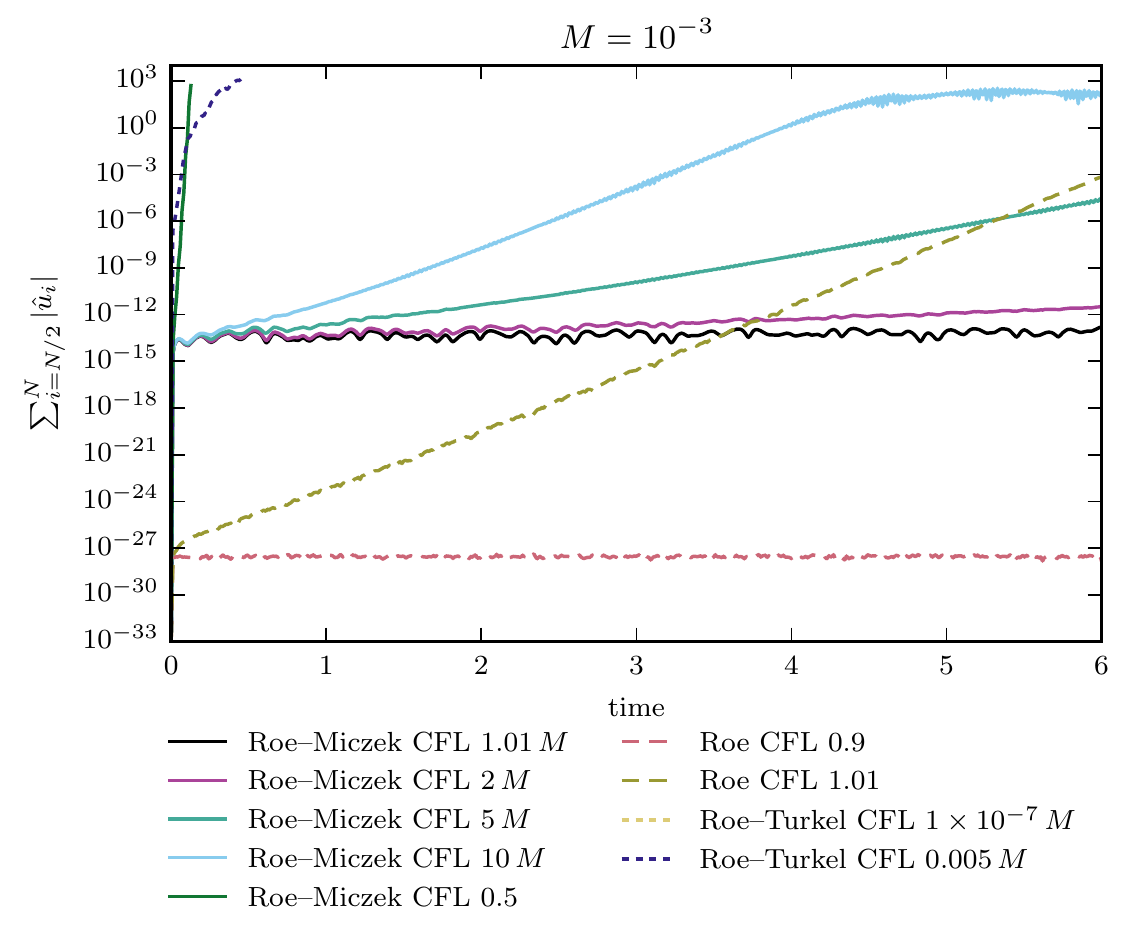}
  \caption{\label{fig:sound-m1e-3}Stability tests of a sound wave at
  Mach number $M=10^{-3}$ using constant reconstruction and forward
  Euler time stepping. The line styles signify the different numerical
  fluxes. For the Roe--Miczek and Roe--Turkel scheme the CFL condition
  was multiplied by $M$ for most tests as indicated.}
\end{figure*}

This more restrictive CFL condition is also observed in the
experiments. As a test setup for CFL stability we use a one-dimensional sound wave.
The initial profile is given by
\begin{align}
  \label{eq:sound1d-initial}
  \rho(x) &= p_0 (1 + M \cos(k x)),\\
  u(x)    &= M c_0 \cos(k x),\\
  p(x)    &= p_0 + \rho_0 c_0^2 M \cos(k x),
\end{align}
with free parameters for background pressure~$p_0$, density~$\rho_0$,
and the corresponding speed of sound~$c_0=\sqrt{\gamma p_0 / \rho_0}$.
The amplitude and thus the Mach number of this sound wave can be adjusted with the
parameter~$M$. The size of the domain is $[0,1)$ with periodic boundary
conditions. We run this setup for a time $t=1$ with explicit
forward Euler time-stepping and constant reconstruction. Explicit integration in time is very inefficient at
low Mach numbers and it can take many million time steps for the
instability to become obvious (i.e.\ visible noise in the velocity
field, negative densities, \ldots). To facilitate the analysis we
compare the growth of the high-frequency Fourier modes, for which we
expect exponential growth in the unstable case. We use
\begin{equation}
  \label{eq:sound1d-mode}
  \sum_{i=N/2}^N |\hat u_{i}|,
\end{equation}
where $\hat u$ is the discrete Fourier transform of $u(x)$ and $N$ is the
number of grid points. We test the growth of this quantity for a range
of $M$ and some values of $\Delta t$ above and below the critical CFL
time step, both for the usual, and the modified Roe scheme.
To emphasize the effect of the numerical flux, we intentionally choose
forward Euler time stepping and constant reconstruction of the interface
values. The results of this experiment at $M=10^{-3}$ are summarized in
Fig.~\ref{fig:sound-m1e-3}. The tests using the standard Roe scheme
confirm that the stability threshold is at CFL~1, as expected. For the
Roe--Miczek scheme we need an additional factor of $M$ in the time step
criterion as it was shown above. The Roe--Turkel scheme was not stable
for any of the tested time steps.

The stability of the discretization of \eqref{eq:euler1}--\eqref{eq:euler3} for the particular case of checkerboard modes ($\beta = \pi$) has equally been confirmed by experiments.

\section{Implicit time discretization}
\label{sect:implicit}

  In the light of the stricter CFL condition of some of the low Mach
  schemes mentioned above, it is natural to turn to implicit time
  discretizations to allow larger time steps. Even for time-explicit schemes with the
  standard CFL condition time steps become prohibitively small: they
  scale with the inverse of the sound speed and therefore are bound to the
  \emph{acoustic} time scale, i.e.\ the sound-crossing time in one
  computational grid cell. In contrast, the
  criterion for selecting the time step in implicit schemes is not
  derived from considerations of stability of the scheme but
  rather from the intended accuracy of the solution. In the low Mach
  case one is usually interested in phenomena that are associated with the fluid
  flow instead of sound waves. In order to accurately resolve the flow,
  the time step should be restricted to the flow crossing time over a
  grid cell -- the \emph{advective} time step criterion. The ratio
  between the acoustic and the advective time steps (and thus the ratio
  of time steps to be taken for bridging the same physical time
  interval) is a function of the Mach number. Even considering the
  increased computational cost of implicit time steps, it is expected
  that there is a certain Mach number below which an implicit scheme is
  more efficient than explicit time discretization. But the threshold
  (and the very feasibility of an implicit time integration) depends on
  the system of equations to be solved and on the efficiency of the
  solution method.

  Implicit time stepping for the Euler equations involves the solution
  of a large nonlinear system of equations. In the three-dimensional
  case the number of equations is $5 N_x N_y N_z$, where is $N_{x,y,z}$
  is the number of grid cells in  $x$, $y$, or $z$ direction. Even for a
  moderately sized problem of $512^3$ cells this already yields about
  $6.7\times 10^8$ nonlinear coupled equations with the same number of unknowns. We use the Newton--Raphson method
  for the solution, which itself requires the solution of a large linear system of
  equations given by the Jacobian of the nonlinear one. Apart from
  the particular implementation approach, the success and efficiency of the
  scheme critically depend on the structure of the system of equations
  to be solved. In particular, a high
  \emph{condition number} of this system would severely impede the
  ability to find solutions efficiently. The particular definition of
  the condition number~$\kappa$, that we use here, is,

  \begin{equation}
    \label{eq:cond-def}
    \kappa = \| A \|_1 \cdot \| A^{-1} \|_1,
  \end{equation}
  where $A$ is the Jacobian matrix of the nonlinear system. We use the
  1-norm as is it computationally less expensive to evaluate compared to
  the 2-norm but still has similar significance for the solution
  efficiency of the linear system.

\begin{figure*}
  \centering
  \includegraphics{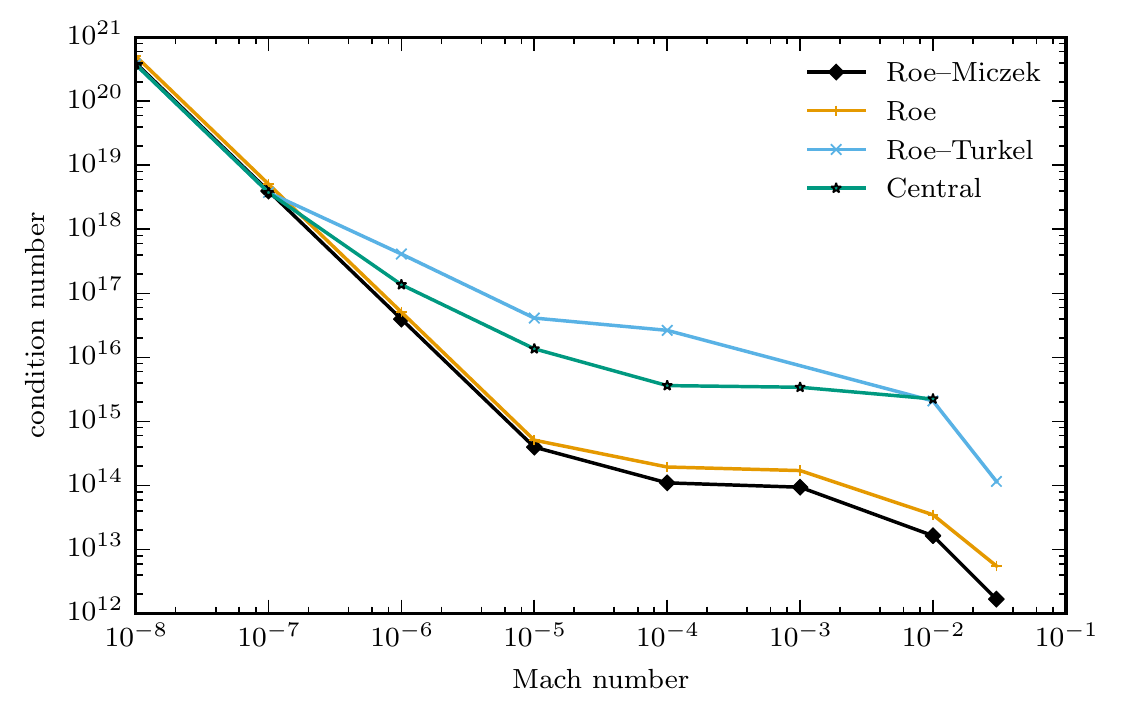}
  \caption{\label{fig:cond-sforce}Mach number dependence of the
  condition number of the unpreconditioned Jacobian matrix that occurs in the solution of
  the nonlinear system. It was obtained on a $16^3$ grid in the
  simulation of stochastically driven turbulence.}
\end{figure*}

  In Fig.~\ref{fig:cond-sforce} we show the effects on the condition number
  for the different modified diffusion matrices presented in this article. In order
  to get a representative condition number for different discretizations
  under realistic conditions we pick a turbulent flow field that was
  produced using stochastic forcing \cite{eswaran1988a,schmidt2006e}.  The first obvious
  feature of the curves is that they are almost identical in the high
  and low Mach number limit. In the regime $M \lesssim 10^{-6}$ this is
  due to the dominating influence of the central flux terms, which all
  the other schemes also include. In the Mach number regime from about
  $10^{-4}$ to $10^{-2}$ the condition number of the Roe-type schemes is
  almost Mach number independent, with the condition number of the
  Roe--Turkel scheme being significantly higher than the other two. For
  practical applications using implicit time stepping this means that
  the Roe--Miczek scheme is as efficient as the standard Roe scheme
  while still giving an accurate result like the Roe--Turkel scheme. Moreover it has been demonstrated \cite{hammer2015a} that the corresponding implementation scales up to $\sim$100,000 cores making it suitable for highly resolved simulations.

\section{Application to unsteady low Mach number flows}
The application examples demonstrated thus far are steady state
problems. The fact that a scheme behaves well at low Mach numbers for
these problems might not be representative of good low Mach behavior in
general. Because of this we apply the Roe--Miczek scheme to a number of
unsteady flow problems at low Mach numbers.

The \emph{Taylor--Green vortex} (TGV) \cite{taylor1937a} is a
large-scale, three-dimensional vortex that decays to smaller vortices,
thereby creating a turbulent flow pattern. Its advantages are that it
can be easily implemented by just setting an initial condition and it
provides a simple estimate for the Reynolds number. In the context of
the Euler equations we can use it to measure numerical viscosity.

We use the initial conditions given in \cite{drikakis2007a},
\begin{equation}
\begin{aligned}
  \rho(t=0) &= \rho_0 = 1.178 \times 10^{-3},\\
u(t=0) &=   u_0 \sin (kx) \cos (ky) \cos (kz),\\
v(t=0) &= - u_0 \cos (kx) \sin (ky) \cos (kz),\\
w(t=0) &= 0,\\
u_0 &= 10^{4},\\
k &= 10^{-2}\\
p(t=0) &= p_0 + \left[{u_0}^2\rho / 16\right]\left[2+\cos \frac{2z}{100}\right]\left[\cos \frac{2x}{100} + \cos \frac{2y}{100}\right],\\
p_0 &= 10^6.
\end{aligned}
\label{eq:tgvortex-initial}
\end{equation}

For comparability with other results in the literature we choose the
value \mbox{$\gamma = 1.4$} in the equation of state (Eq.~\eqref{eq:eos}).
The maximum Mach number of this setup is $M_\text{max} = u_0 / \sqrt{\gamma
p_0 / \rho_0} \approx 0.29$. We can easily scale this setup to lower
Mach numbers by multiplying $u_0$ with the appropriate factor.

To be able to compare simulations at different Mach numbers, we
scale certain quantities (denoted by $^*$). The relations are
\begin{equation}
  t^* = k u_0 t, \quad K^* = K / {u_0}^2, \quad \Omega^* = \Omega / (k u_0)^2.
  \label{eq:tgvortex-nondim}
\end{equation}

To test the effect of the diffusion matrix modification we
calculate the numerical Reynolds number of the usual Roe scheme and the modified one at different resolutions and Mach numbers. The
Reynolds number is purely numerical as we do not include any explicit
viscosity terms.  We use the expression given in \cite{taylor1937a},
\begin{equation}
  \frac{dK^*}{dt^*} = - \frac{\Omega^*}{\Reyn},
  \label{eq:tgvortex-Re}
\end{equation}
with the mean of kinetic energy~$K$ and the mean enstrophy~$\Omega$,
\begin{equation}
  K = \frac{1}{2} \langle |\fvec{v}|^2 \rangle, \quad \Omega =
  \frac{1}{2} \langle |\nabla \times \fvec{v}|^2 \rangle.
  \label{eq:tgvortex-K-Omega}
\end{equation}
The operation $\langle \cdot \rangle$ is a volumetric average. Both
quantities can be independently computed from the velocity
field~$\fvec{v}$.

\begin{figure*}
  \centering
  \begin{tabular}{ccc}
    \includegraphics[width=0.3\textwidth]{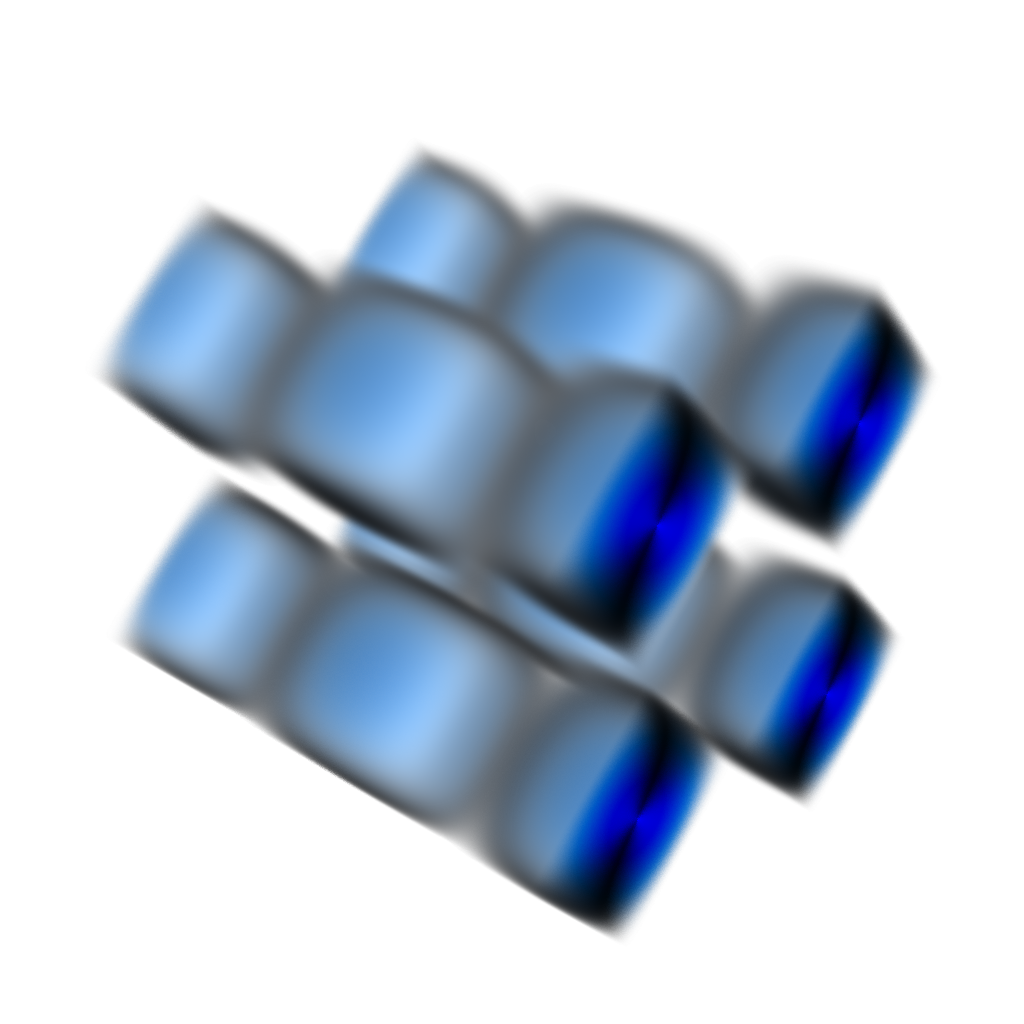}&
    \includegraphics[width=0.3\textwidth]{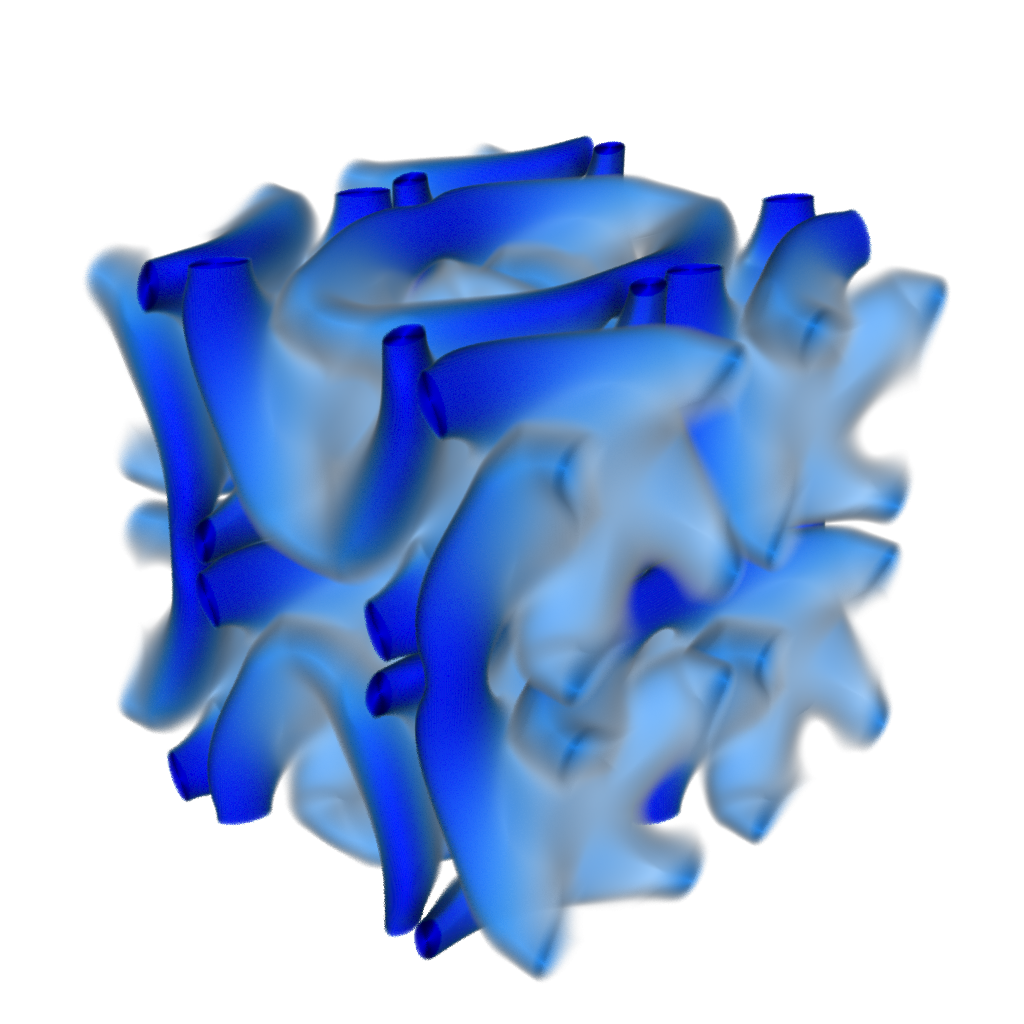}&
    \includegraphics[width=0.3\textwidth]{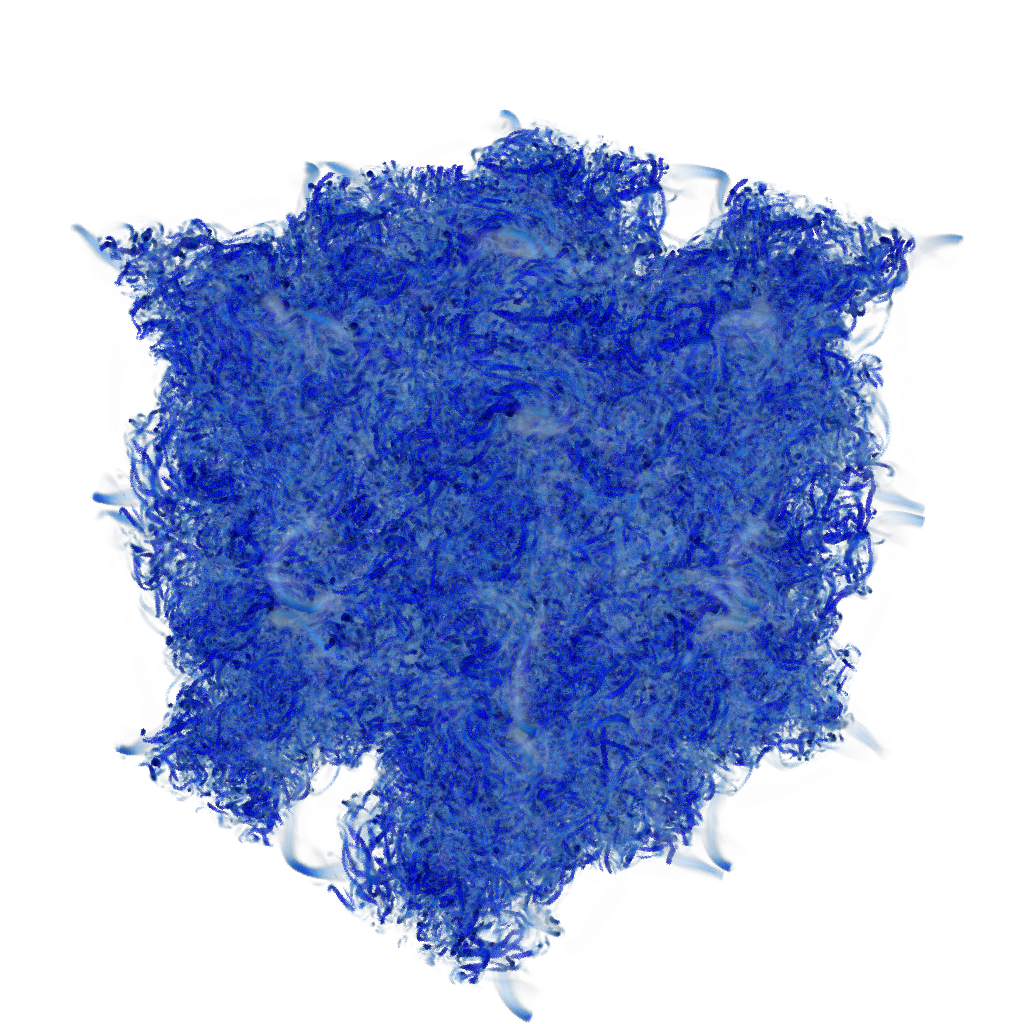}\\
    $t^* = \num{0}$&
    $t^* = \num{0.96}$&
    $t^* = \num{3.05}$\\
    scale $10$&
    scale $10$&
    scale $1000$\\
    \includegraphics[width=0.3\textwidth]{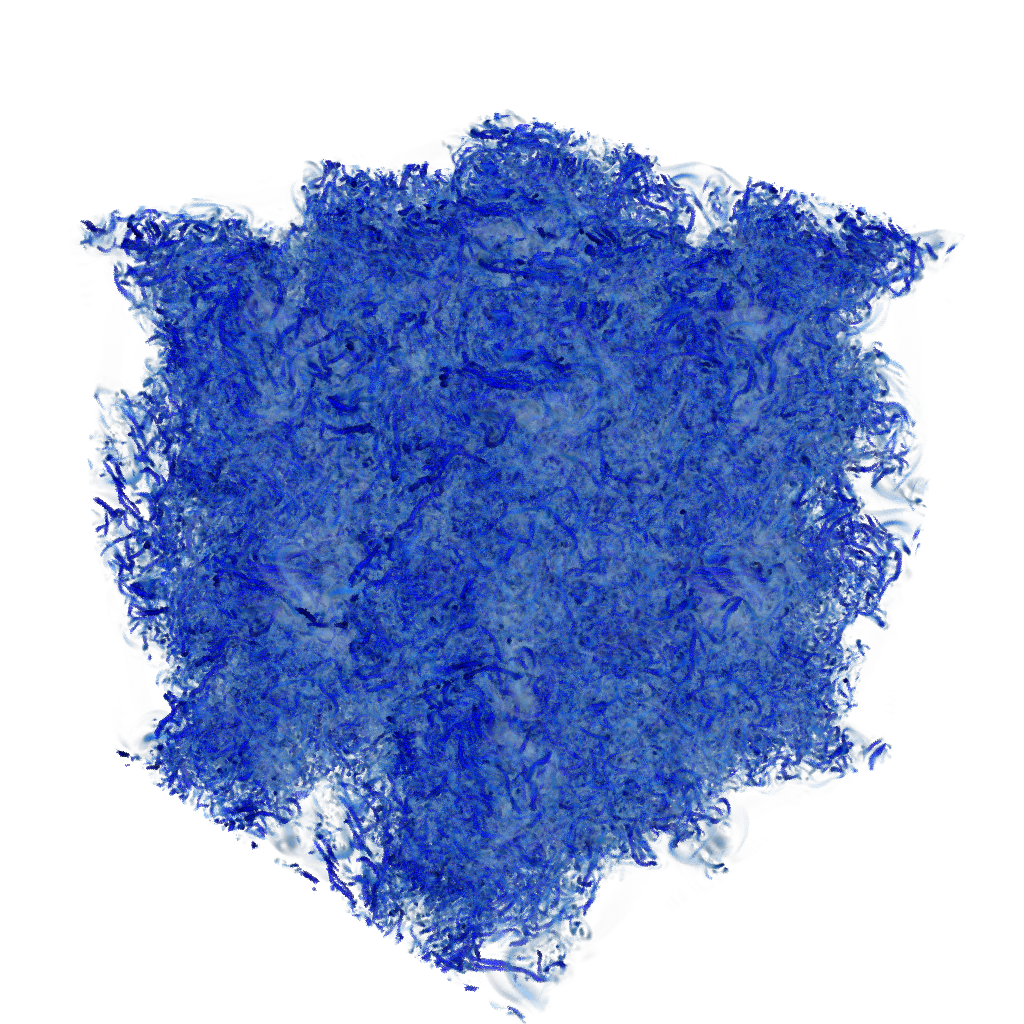}&
    \includegraphics[width=0.3\textwidth]{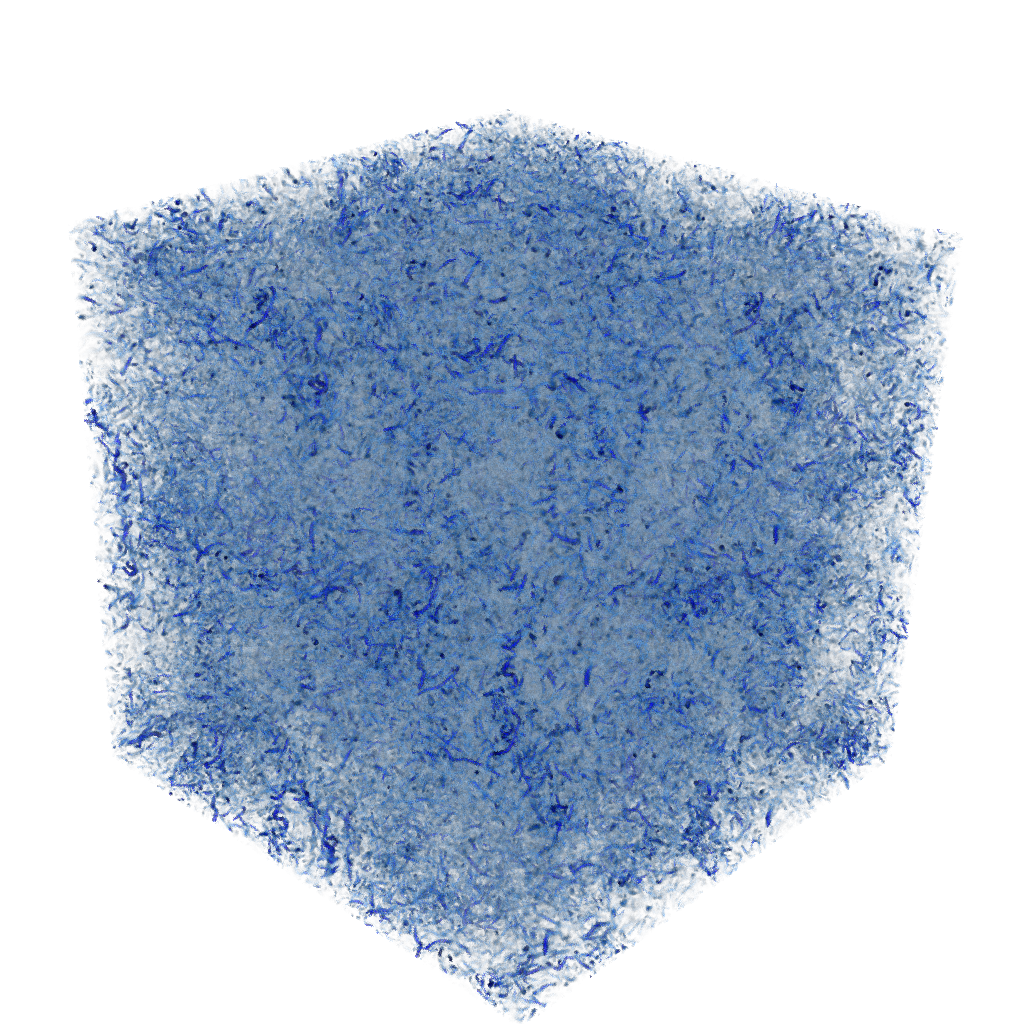}&
    \includegraphics[width=0.3\textwidth]{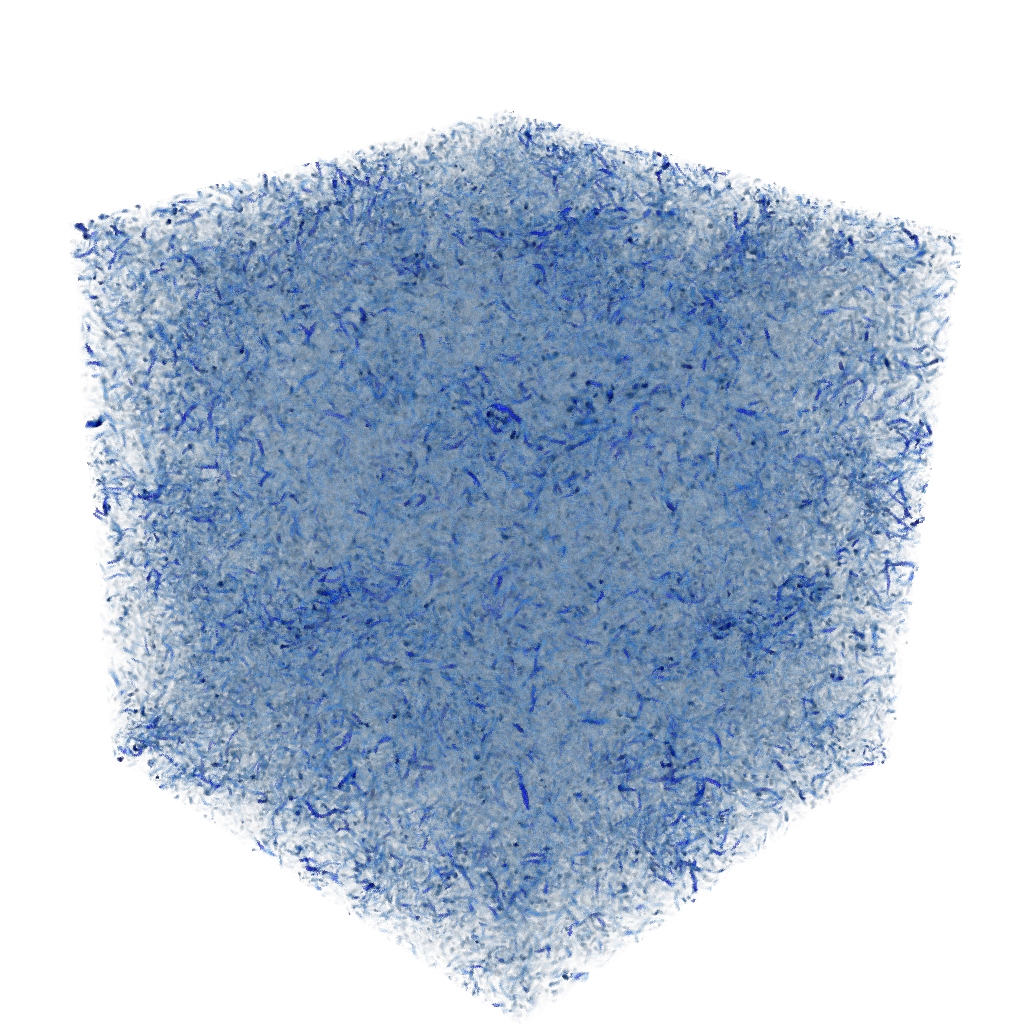}\\
    $t^* = \num{3.34}$&
    $t^* = \num{5.01}$&
    $t^* = \num{8.70}$\\
    scale 1000&
    scale 1000&
    scale 250
  \end{tabular}
  \caption{\label{fig:tgv-l2}Temporal evolution of vortex cores in the Taylor--Green
  vortex visualized using a criterion from \cite{jeong1995a}. The Mach
  number of the initial configuration was set to $10^{-2}$. The
  simulation was run at a resolution of $512^3$ grid cells using the Roe--Miczek
  scheme and implicit ESDIRK34 time stepping. The magnitude of the color
  scale was adjusted as stated below each panel.}
\end{figure*}

Figure~\ref{fig:tgv-l2} visualizes the decay of the large scale
vortices to smaller scales and their eventual dissipation using a
criterion from \cite{jeong1995a}.

\begin{figure*}
  \includegraphics{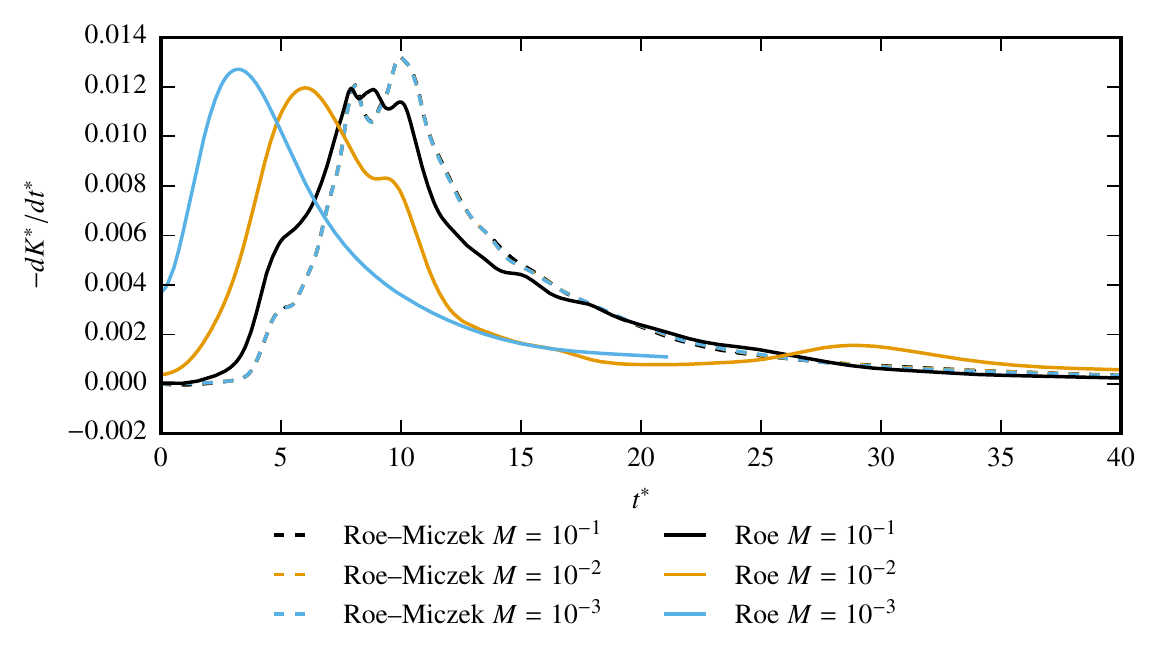}
  \caption{\label{fig:tgv-128}Taylor--Green vortex
  simulated at a fixed resolution of $128^3$ grid cells. The solid lines
  were computed using the Roe flux, the dashed lines using the
  modified flux \cite{miczek2015a}. The colors signify different
  initial Mach numbers.}
\end{figure*}

The test in Fig.~\ref{fig:tgv-128} shows kinetic energy dissipation rate
of the vortex at different Mach numbers computed with a fixed grid size
of $128^3$ cells. It is expected that this rate is independent of Mach
number in the nondimensional variables of
Eq.~\eqref{eq:tgvortex-nondim}. This is very well fulfilled for the
modified Roe scheme from \cite{miczek2015a} but not for the original Roe
scheme.

\section{Tests in the High-Mach number regime}
\label{sec:high-mach}
\begin{figure*}
  \includegraphics{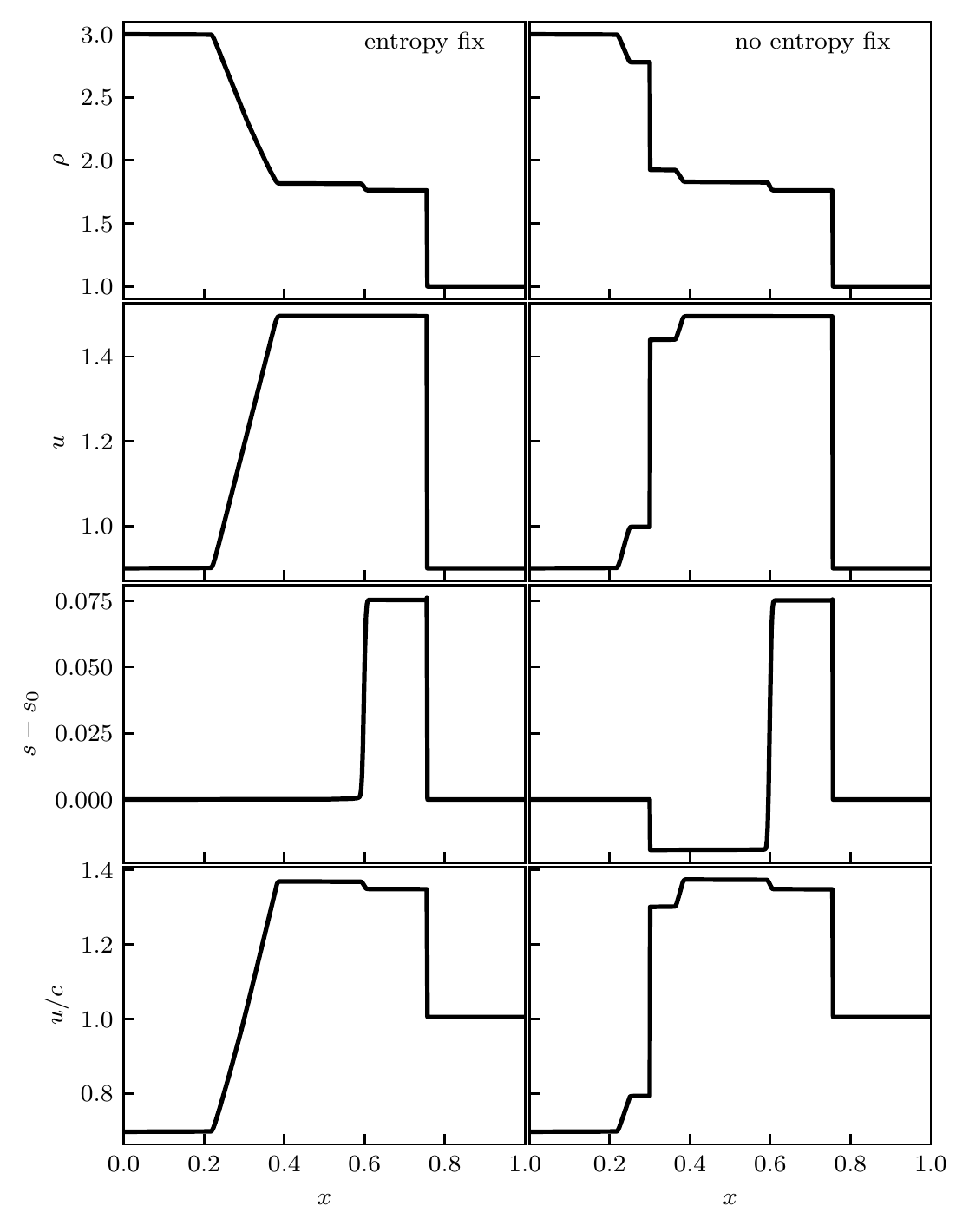}
  \caption{\label{fig:entropyfix}This shows density (\emph{top row}),
  velocity (\emph{second row}), deviation of entropy from its initial value
  (\emph{third row}), and Mach number (\emph{bottom row}) of a transonic shock
  tube defined by the initial conditions in Eq.~\eqref{eq:shock-tube} at time
  $t=0.2$. The right column was computed using the Roe--Miczek scheme without
  any entropy fix. The left column used the same scheme with the entropy fix
  suggested in \cite{harten1983b}.}
\end{figure*}

To qualify as a scheme for all Mach numbers it has to be
ensured that it behaves correctly in the presence of shocks. The
modifying matrix, Eq.~\eqref{eq:pc_miczek_sym}, continuously
approaches the identity matrix as the Mach number increases to 1. The
newly proposed scheme is identical to the original Roe scheme for $M_\text{loc} \geq 1$. 
The Roe scheme, however, is known to violate the entropy condition in transonic rarefactions, which is cured by introducing an entropy fix (see \cite{pelanti2001a} for an overview). Here we test if this treatment interferes with our low Mach number modifications. For the present study we
chose to implement the fix suggested by \cite{harten1983b}. We use a shock
tube with isentropic initial conditions, which makes the entropy problem and its solution obvious. The initial setup in
primitive variables is
\begin{equation}
  \label{eq:shock-tube}
  (\rho,u,p)(x,0) = \begin{cases}
    (3,0.9,3)            &\text{for $x < 0.5$},\\
    (1,0.9,3^{1-\gamma}) &\text{otherwise}.
  \end{cases}
\end{equation}
The value $\gamma$ is the adiabatic index in the ideal gas equation of state \eqref{eq:eos},
i.e. the internal energy is $e = p / (\gamma -1)$. Figure~\ref{fig:entropyfix} shows the
solution at time $t=0.2$ for the Roe--Miczek scheme with and without the
entropy fix. Implementing an entropy fix shows no interference with the introduced modifications.

\section{Conclusions}

We have presented a finite volume solver for the compressible Euler equations (\cite{miczek2015a}) with the following properties:
\begin{enumerate}[(i)]
\item it is preserving the asymptotics of the low Mach limit
\item it ensures the kinetic energy of the flow is not excessively dissipated by the numerics near the incompressible regime
\item it is a linearly stable scheme
\item the inversion of the large system of equations arising from implicit time discretization has a good condition number.
\end{enumerate}
Requirement (\ref{req:2}) reflects well the issue of growing numerical diffusion in the low Mach limit observed for shock-capturing schemes, and might be easier to check than (\ref{req:1}).

We have shown numerical simulations with this scheme that work for Mach numbers as low as $10^{-10}$.  The method can also be applied to flows with high Mach number, where it continuously approaches the Roe solver. For transonic flows our modification does not interfere with the entropy fix. In this sense the method of Miczek et\,al. \cite{miczek2015a} constitutes an all-Mach number scheme.

We conclude that when devising a low Mach number scheme for the Euler equations it is important that the discrete equations reflect the nature of the limit system for finite discretizations already. It thus is not sufficient to ensure that asymptotically the discrete equations become consistent discretizations of the limit system. We show that in addition one needs to satisfy the requirements \ref{req:2}, \ref{req:3}, \ref{req:4} for such a scheme to be feasible in practice.

In upcoming work we plan to extend these properties to the Euler equations with a gravitational source term, where one needs to ensure in addition that the scheme is well-balanced.

\section*{Appendix}

The limit of low Mach numbers in the context of the Euler equations \eqref{eq:euler1}--\eqref{eq:euler3} is best explored by introducing a family of solutions, parametrized by a real dimensionless number $M > 0$, $M \to 0$. Writing $f \in \mathcal O(M^p)$ means that the leading order of the expansion of $f$ in powers of $M$ is $M^p$, i.e. $f = M^p (f^{(0)} + f^{(1)} M + \ldots )$, with the functions $f^{(i)}$ not depending on $M$.
\\Only two requirements are needed to derive the rescaled Euler equations:
\begin{itemize}
  \item The local Mach number $\displaystyle M_\text{loc}(\fvec x, t) := |\fvec v(\fvec x,t)| / \sqrt{\frac{\gamma p(\fvec x, t)}{\rho(\fvec x,t)}}$ shall be scaling uniformly with $M$ as $M \to 0$: $M_\text{loc} \in \mathcal{O}(M)$.
  \item Every member of the family shall fulfill the same equation of state
\begin{equation}
  \label{eq:eos}
  E = \frac{ p}{\gamma-1} + \frac{1}{2}  \rho | \fvec v|^2.
\end{equation}
\end{itemize}
The most general asymptotic scalings in this case are
\begin{align} 
  \fvec x &= M^{\mathfrak a} \tilde {\fvec x}, & t &=  M^{\mathfrak b} \tilde t,\\
 \rho(\fvec x, t) &= M^{\mathfrak c + 2 - 2 \mathfrak d}\tilde \rho(\fvec {\tilde x}, \tilde t), & \fvec v(\fvec x, t) &= M^{\mathfrak d} \tilde {\fvec v}(\fvec{\tilde x}, \tilde t), \\ 
 E(\fvec x, t) &= M^{\mathfrak c} \tilde E(\fvec {\tilde x}, \tilde t), & p(\fvec x, t) &= M^{\mathfrak c} \tilde p(\fvec {\tilde x}, \tilde t),
\end{align}
with $\mathfrak{a, b, c, d}$ arbitrary numbers, as can be found from a direct computation. It is understood that quantities with a tilde are $\mathcal{O}(1)$ when expanded as power series in $M$. An example of such a family of solutions is given by the Gresho vortex setup in \eqref{eq:gresho1}--\eqref{eq:gresho2}.

Furthermore every member of the family shall fulfill the Euler equations. Inserting the above scalings yields a system of equations that is fulfilled by quantities with a tilde. These equations shall be called \emph{rescaled}, and cannot be the Euler equations again, because the Mach number changes. They are found to be
\begin{equation}
 \tilde E = \frac{\tilde p}{\gamma-1} + \frac{1}{2} M^2 \tilde\rho |\tilde {\fvec v}|^2
\end{equation}
and
\begin{align}
 &M^{\mathfrak a-\mathfrak d-\mathfrak b} \del_t \tilde \rho + \nabla \cdot (\tilde \rho \tilde {\fvec v}) = 0,\\
 &M^{\mathfrak a-\mathfrak d-\mathfrak b} \del_t (\tilde \rho \tilde {\fvec v}) + \nabla \cdot \left(\tilde \rho \tilde {\fvec v} \otimes \tilde {\fvec v} + \frac{\tilde p}{M^2}\cdot \id\right ) = 0,\\
 &M^{\mathfrak a-\mathfrak d-\mathfrak b} \del_t \tilde E + \nabla \cdot (\tilde {\fvec v}(\tilde E+\tilde p)) = 0.
\end{align}
Observe the fact that the kinetic energy obtains an additional factor of $M^2$ in the equation of state.

The factor in front of the time derivatives is related to the dimensionless Strouhal number 
\begin{align}
 \text{\textit{Str}}_\text{loc} &= \frac{x}{|\fvec v| t} = \frac{M^{\mathfrak a} \tilde x}{M^{\mathfrak d} \tilde {\fvec v} \cdot M^{\mathfrak b} \tilde t}.
\end{align}
This factor is \emph{not} identical to the Strouhal number, but is just its asymptotic $M$-scaling. As an additional condition on the family of solutions one is tempted to insist on $\text{\textit{Str}} \in \mathcal{O}(1)$, i.e. $\mathfrak a - \mathfrak d - \mathfrak b = 0$. This corresponds to adapting the time scales to the speed of the fluid (and not to sound wave crossing times). 
\\Different ways of decreasing the Mach number (e.g. by decreasing the value of the velocity, or increasing the sound speed instead, or a combination of both) are equivalent and should result in the same rescaled equations. This explains why the precise value of $\mathfrak{a, b, c, d}$ does not matter for the form of the rescaled equations. Only these equations will be considered in what follows and we drop the tilde.

\begin{acknowledgements}
  We thank Philipp Birken for stimulating discussions. WB gratefully acknowledges support from the German National Academic
  Foundation. The work of FKR and PVFE was supported by the
  Klaus Tschira Foundation. CK acknowledges support of the DFG priority program SPPEXA. The authors gratefully acknowledge the
  Gauss Centre for Supercomputing (GCS) for providing computing time
  through the John von Neumann Institute for Computing (NIC) on the
  GCS share of the supercomputer JUQUEEN~\cite{stephan2015a} at J{\"u}lich Supercomputing
  Centre (JSC). GCS is the alliance of the three national supercomputing
  centres HLRS (Universit{\"a}t Stuttgart), JSC (Forschungszentrum J{\"u}lich),
  and LRZ (Bayerische Akademie der Wissenschaften), funded by the German
  Federal Ministry of Education and Research (BMBF) and the German State
  Ministries for Research of Baden-W{\"u}rttemberg (MWK), Bayern (StMWFK) and
  Nordrhein-Westfalen (MIWF). 
\end{acknowledgements}

\bibliographystyle{spmpsci}      

\end{document}